# UNIRATIONALITY OF DEL PEZZO SURFACES OF DEGREE TWO OVER FINITE FIELDS (EXTENDED VERSION)

DINO FESTI, RONALD VAN LUIJK

ABSTRACT. We prove that every del Pezzo surface of degree two over a finite field is unirational, building on the work of Manin and an extension by Salgado, Testa, and Várilly-Alvarado, who had proved this for all but three surfaces. Over general fields, we state several sufficient conditions for a del Pezzo surface of degree two to be unirational.

## 1. INTRODUCTION

We say that a variety $X$ of dimension $n$ over a field $k$ is *unirational* if there exists a dominant rational map $\mathbb{P}^n \dashrightarrow X$, defined over $k$. A *del Pezzo surface* is a smooth, projective, geometrically integral variety $X$ of which the anticanonical divisor $-K_X$ is ample. We define the *degree* of a del Pezzo surface $X$ as the self intersection number of $K_X$, that is, $\deg X = K_X^2$. If $k$ is an algebraically closed field, then every del Pezzo surface of degree $d$ over $k$ is isomorphic to $\mathbb{P}^1 \times \mathbb{P}^1$ (with $d = 8$), or to $\mathbb{P}^2$ blown up in $9 - d$ points in general position.

Over arbitrary fields, the situation is more complicated. Works of B. Segre, Yu. Manin, J. Kollár, M. Pieropan, and A. Knecht prove that every del Pezzo surface of degree $d \geq 3$, defined over any field $k$, is unirational, provided that the set $X(k)$ of rational points is non-empty. For references, see [Seg43, Seg51] for $k = \mathbb{Q}$ and $d = 3$, see [Man86, Theorem 29.4 and 30.1] for $d \geq 3$ with the extra assumption for $d \in \{3, 4\}$ that $k$ has enough elements. See [Kol02, Theorem 1.1] for $d = 3$ and a general ground field. The earliest reference we could find for $d = 4$ and a general ground field is [Pie12, Proposition 5.19]. Independently, for $d = 4$, [Kne13, Theorem 2.1] covers all finite fields.

Since all del Pezzo surfaces over finite fields have a rational point (see [Man86, Corollary 27.1.1]), this implies that every del Pezzo surface of degree at least 3 over a finite field is unirational.

Building on work by Manin (see [Man86, Theorem 29.4]), C. Salgado, D. Testa, and A. Várilly-Alvarado prove that all del Pezzo surfaces of degree 2 over a finite field are unirational as well, except possibly for three isomorphism classes of surfaces (see [STVA14, Theorem 1]). In this paper, we show that these remaining three cases are also unirational, thus proving our first main theorem.

**Theorem 1.1.** *Every del Pezzo surface of degree 2 over a finite field is unirational.*

In Section 2 we will recall the basics about del Pezzo surfaces of degree 2, including the fact that the linear system associated to the anti-canonical divisor induces a finite morphism to $\mathbb{P}^2$ of degree 2. This allows us to state the second main theorem.

**Theorem 1.2.** *Suppose $k$ is a field of characteristic not equal to 2, and let $\overline{k}$ be an algebraic closure of $k$. Let $X$ be a del Pezzo surface of degree 2 over $k$. Let $B \subset \mathbb{P}^2$ be the branch locus of the anti-canonical morphism $\pi \colon X \to \mathbb{P}^2$. Let $C \subset \mathbb{P}^2$ be a projective curve that is birationally equivalent with $\mathbb{P}^1$ over $k$. Assume that all singular points of $C$ that are contained in $B$ are ordinary singular points. Then the following statements hold.*

  (1) *Suppose that there is a point $P \in X(k)$ such that we have $Q = \pi(P) \in C - B$. Suppose that $B$ contains no singular points of $C$ and that all intersection points of $B$ and $C$ have even intersection multiplicity. Then the surface $X$ is unirational.*
  (2) *Suppose that one of the following two conditions hold.*

(a) There is a point $Q \in C(k) \cap B(k)$ that is a double or a triple point of $C$. The curve $B$ contains no other singular points of $C$, and all intersection points of $B$ and $C$ have even intersection multiplicity.

(b) There exist two distinct points $Q_1, Q_2 \in C(\overline{k}) \cap B(\overline{k})$ such that $B$ and $C$ intersect with odd multiplicity at $Q_1$ and $Q_2$ and with even intersection multiplicity at all other intersection points. Furthermore, the points $Q_1$ and $Q_2$ are smooth points or double points on the curve $C$, and $B$ contains no other singular points of $C$.

Then there exists a field extension $\ell$ of $k$ of degree at most $2$ for which the preimage $\pi^{-1}(C_\ell)$ is birationally equivalent with $\mathbb{P}^1_\ell$; for each such field $\ell$, the surface $X_\ell$ is unirational.

**Corollary 1.3.** *Suppose $k$ is a field of characteristic not equal to $2$. Let $X$ be a del Pezzo surface of degree $2$ over $k$. Assume that $X$ has a $k$-rational point, say $P$. Let $C \subset \mathbb{P}^2$ be a geometrically integral curve over $k$ of degree $d \geq 2$ and suppose that $\pi(P)$ is a point of multiplicity $d-1$ on $C$. Suppose, moreover, that $C$ intersects the branch locus $B$ of the anti-canonical morphism $\pi \colon X \to \mathbb{P}^2$ with even multiplicity everywhere. Then the following statements hold.*

(1) *If $\pi(P)$ is not contained in $B$, then $X$ is unirational.*
(2) *If $\pi(P)$ is contained in $B$, and it is an ordinary singular point on $C$ and we have $d \in \{3, 4\}$, then there exists a field extension $\ell$ of $k$ of degree at most $2$ for which the preimage $\pi^{-1}(C_\ell)$ is birationally equivalent with $\mathbb{P}^1_\ell$; for each such field $\ell$, the surface $X_\ell$ is unirational.*

After a quick summary of some basic results on del Pezzo surfaces of degree 2 in the next section, we will present the three difficult surfaces and prove Theorem 1.1 in Section 3. The main tool is Lemma 3.2 (that is, [STVA14, Theorem 17]), which states that it suffices to construct a rational curve on each of the three del Pezzo surfaces.

In Section 4 we give an alternative description of Manin's construction of a rational curve. If a point $P$ on a del Pezzo surface $X$ does not lie on four of the 56 exceptional curves, and not on the ramification locus of the anti-canonical morphism, then Manin's construction, extended by C. Salgado, D. Testa, and A. Várilly-Alvarado, yields a rational curve that satisfies the assumptions of case (1) of Corollary 1.3 with the degree $d$ being such that there are $4-d$ exceptional curves through $P$ (cf. Proposition 4.14, Remark 4.15, and Example 5.14).

The three difficult surfaces do not contain such a point. The proofs of unirationality of these three cases use a rational curve that is an example of case (2) of Corollary 1.3 instead (cf. Remark 3.4 and Example 5.16). Here we benefit from the fact that if $k$ is a finite field, then any curve that becomes birationally equivalent with $\mathbb{P}^1$ over an extension of $k$, already is birationally equivalent with $\mathbb{P}^1$ over $k$ itself. For two of the three cases, the rational curve we use has degree 4. For the last case, the curve we use has degree 3, but there also exist quartic curves satisfying the hypotheses of case (2) of Corollary 1.3. This raises the following question (cf. Question 6.6, Remark 6.8, and Example 6.9), which together with case (2) of Corollary 1.3 could help proving unirationality of del Pezzo surfaces of degree 2 over any field of characteristic not equal to 2.

**Question 1.4.** *Let $d \in \{3, 4\}$ be an integer. Let $X$ be a del Pezzo surface of degree two over a field of characteristic not equal to $2$, and let $P \in X(k)$ be a point on the ramification locus of the anti-canonical map $\pi \colon X \to \mathbb{P}^2$. Does there exist a geometrically integral curve of degree $d$ in $\mathbb{P}^2$ over $k$ that has an ordinary singular point of multiplicity $d-1$ at $\pi(P)$, and that intersects the branch locus of $\pi$ with even multiplicity everywhere?*

For some $d$, $X$, and $P$, the answer to this question is negative (see Example 6.9), but in all cases we know of (all over finite fields), there do exist degenerate curves of degree $d$ with a point of multiplicity at least $d-1$ at $\pi(P)$. Hence, it may be true that the answer to Question 1.4 is positive for $X$ and $P$ general enough.

In line with case (1) of Corollary 1.3, we can ask, in fact for any integer $d \geq 1$, an analogous question for points $P$ that do not lie on the ramification locus, where we do not require the singular point to be ordinary. In this case, if $P$ lies on $r \leq 3$ exceptional curves, then Manin's construction shows that the answer is positive for degree $d = 4 - r$. Therefore, this analogous question is especially interesting when $P$ lies on four exceptional curves (cf. Remark 5.12 and Example 5.15).

In Section 5 we prove Theorem 1.2 and a generalisation, Corollary 5.6, as well as Corollary 1.3. We do this by showing that, under the hypotheses of Theorem 1.2, the pull-back of the curve $C$ to $X$ contains a rational component, so that we may again use Lemma 3.2.

In Section 6 we discuss how to search for curves satisfying the assumptions of Theorem 1.2 and in particular of Corollary 1.3.

The authors would like to thank David Holmes, Ariyan Javanpeykar, Bjorn Poonen, Marco Streng, Damiano Testa, and Anthony Várilly-Alvarado for useful conversations.

## 2. Del Pezzo surfaces of degree two

The statements in this section are well known and we will use them freely. Let $X$ be a del Pezzo surface of degree 2 over a field $k$ with canonical divisor $K_X$. The Riemann-Roch spaces $\mathcal{L}(-K_X)$ and $\mathcal{L}(-2K_X)$ have dimension 3 and 7, respectively. Let $x, y, z$ be generators of $\mathcal{L}(-K_X)$ and choose an element $w \in \mathcal{L}(-2K_X)$ that is not contained in the image of the natural map $\mathrm{Sym}^2 \mathcal{L}(-K_X) \to \mathcal{L}(-2K_X)$. Then $X$ embeds into the weighted projective space $\mathbb{P} = \mathbb{P}(1,1,1,2)$ with coordinates $x, y, z$, and $w$. We will identify $X$ with its image in $\mathbb{P}$, which is a smooth surface of degree 4. Conversely, every smooth surface of degree 4 in $\mathbb{P}$ is a del Pezzo surface of degree 2. There are homogeneous polynomials $f, g \in k[x, y, z]$ of degrees 2 and 4, respectively, such that $X \subset \mathbb{P}$ is given by

$$w^2 + fw = g. \tag{1}$$

If the characteristic of $k$ is not 2, then after completing the square on the left-hand side, we may assume $f = 0$. For more details and proofs of these facts, see [Kol96, Section III.3, Theorem III.3.5] and [Man86, Section IV.24].

The restriction to $X$ of the 2-uple embedding $\mathbb{P} \to \mathbb{P}^6$ corresponds to the complete linear system $|-2K_X|$. Every hyperplane section of $X \subset \mathbb{P}$ is linearly equivalent with $-K_X$. The projection $\mathbb{P} \dashrightarrow \mathbb{P}^2$ onto the first three coordinates restricts to a finite, separable morphism $\pi_X \colon X \to \mathbb{P}^2$ of degree 2, which corresponds to the complete linear system $|-K_X|$.

The morphism $\pi_X$ is ramified above the branch locus $B_X \subset \mathbb{P}^2$ given by $f^2 + 4g = 0$. If the characteristic of $k$ is not 2, then $B_X$ is a smooth curve. We denote the ramification locus $\pi^{-1}(B_X)$ of $\pi_X$ by $R_X$. As for every double cover, the morphism $\pi_X$ induces an involution $\iota_X \colon X \to X$ that sends a point $P \in X$ to the unique second point in the fiber $\pi_X^{-1}(\pi_X(P))$, or to $P$ itself if $\pi_X$ is ramified at $P$. If $X$ is clear from the context, then we sometimes leave out the subscript and write $\pi, \iota, B$, and $R$ for $\pi_X, \iota_X, B_X$, and $R_X$, respectively.

As in [STVA14], we define a *bitangent* to $B$ to be a line $\ell \subset \mathbb{P}^2$ for which $\pi^{-1}(\ell)$ is reducible. If the characteristic of $k$ is not 2, then this is equivalent to a more intuitive definition: a bitangent line is a line in $\mathbb{P}^2$ whose intersection points with $B$ all have even intersection multiplicity. For each bitangent $\ell$ to $B$, the two components of $\pi^{-1}(\ell)$ are exceptional curves, say $E$ and $\iota_X(E)$, whose sum $E + \iota_X(E)$ is linearly equivalent with $-K_X$ (see [STVA14, Lemma 4]). There are exactly 28 bitangent lines to $B$, and the irreducible components above them are exactly the 56 exceptional curves on $X$. Any point of $X$ lies on at most four exceptional curves (see [TVAV09, proof of Lemma 4.1]); a *generalised Eckardt point* is a point on $X$ that lies on exactly four exceptional curves.

## 3. Proof of the first main theorem

Set $k_1 = k_2 = \mathbb{F}_3$ and $k_3 = \mathbb{F}_9$. Let $\gamma \in k_3$ denote an element satisfying $\gamma^2 = \gamma + 1$. Note that $\gamma$ is not a square in $k_3$. For $i \in \{1, 2, 3\}$, we define the surface $X_i$ in $\mathbb{P} = \mathbb{P}(1,1,1,2)$ with coordinates $x, y, z, w$ over $k_i$ by

$$\begin{aligned} X_1: \ -w^2 &= (x^2+y^2)^2 + y^3 z - yz^3, \\ X_2: \ -w^2 &= x^4 + y^3 z - yz^3, \\ X_3: \ \gamma w^2 &= x^4 + y^4 + z^4. \end{aligned}$$

These surfaces are smooth, so they are del Pezzo surfaces of degree 2. C. Salgado, D. Testa, and A. Várilly-Alvarado proved the following result.

**Theorem 3.1.** *Let $X$ be a del Pezzo surface of degree 2 over a finite field. If $X$ is not isomorphic to $X_1, X_2,$ and $X_3$, then $X$ is unirational.*

*Proof.* See [STVA14, Theorem 1]. □

We will use the following lemma to prove the complementary statement, namely that $X_1, X_2,$ and $X_3$ are unirational as well.

**Lemma 3.2.** *Let $X$ be a del Pezzo surface of degree 2 over a field $k$. Suppose that $\rho\colon \mathbb{P}^1 \to X$ is a nonconstant morphism; if the characteristic of $k$ is 2 and the image of $\rho$ is contained in the ramification divisor $R_X$, then assume also that the field $k$ is perfect. Then $X$ is unirational.*

*Proof.* See [STVA14, Theorem 17]. □

For $i \in \{1, 2, 3\}$, we define a morphism $\rho_i\colon \mathbb{P}^1 \to X_i$ by extending the map $\mathbb{A}^1(t) \to X_i$ given by
$$t \mapsto (x_i(t) : y_i(t) : z_i(t) : w_i(t)),$$
where

$$\begin{aligned}
x_1(t) &= t^2(t^2-1), & x_2(t) &= t(t^2+1)(t^4-1), & x_3(t) &= (t^4+1)(t^2-\gamma^3), \\
y_1(t) &= t^2(t^2-1)^2, & y_2(t) &= -t^4, & y_3(t) &= (t^4-1)(t^2+\gamma^3), \\
z_1(t) &= t^8-t^2+1, & z_2(t) &= t^8+1, & z_3(t) &= (t^4+\gamma^2)(t^2-\gamma), \\
w_1(t) &= t(t^2-1)(t^4+1)(t^8+1), & w_2(t) &= t^2(t^2+1)(t^{10}-1), & w_3(t) &= \gamma^2 t(t^8-1)(t^2+\gamma).
\end{aligned}$$

It is easy to check for each $i$ that the morphism $\rho_i$ is well defined, that is, the polynomials $x_i, y_i, z_i,$ and $w_i$ satisfy the equation of $X_i$, and that $\rho_i$ is non-constant.

**Theorem 3.3.** *The del Pezzo surfaces $X_1, X_2,$ and $X_3$ are unirational.*

*Proof.* By Lemma 3.2, the existence of $\rho_1, \rho_2,$ and $\rho_3$ implies that $X_1, X_2,$ and $X_3$ are unirational. □

*Proof of Theorem 1.1.* This follows from Theorems 3.1 and 3.3. □

*Remark* 3.4. Take any $i \in \{1, 2, 3\}$. Set $A_i = \rho_i(\mathbb{P}^1)$ and $C_i = \pi_i(A_i)$, where $\pi_i = \pi_{X_i}\colon X_i \to \mathbb{P}^2$ is as described in the previous section. By Remark 2 of [STVA14], the surface $X_i$ is minimal, and the Picard group $\operatorname{Pic} X_i$ is generated by the class of the anticanonical divisor $-K_{X_i}$. The same remark states that the linear system $|-nK_{X_i}|$ does not contain a geometrically integral curve of geometric genus zero for $n \leq 3$ if $i \in \{1, 2\}$, nor for $n \leq 2$ if $i = 3$. For $i \in \{1, 2\}$, the curve $A_i$ has degree 8, so it is contained in the linear system $|-4K_{X_i}|$. The curve $A_3$ has degree 6, so it is contained in the linear system $|-3K_{X_i}|$. This means that the curve $A_i$ has minimal degree among all rational curves on $X_i$. The restriction of $\pi_i$ to $A_i$ is a double cover $A_i \to C_i$. The curve $C_i \subset \mathbb{P}^2$ has degree 4 for $i \in \{1, 2\}$ and degree 3 for $i = 3$, and $C_i$ is given by the vanishing of $h_i$, with

$$\begin{aligned}
h_1 &= x^4 + xy^3 + y^4 - x^2yz - xy^2z, \\
h_2 &= x^4 - x^2y^2 - y^4 + x^2yz + yz^3, \\
h_3 &= x^2y + xy^2 + x^2z - xyz + y^2z - xz^2 - yz^2 - z^3.
\end{aligned}$$

For $i \in \{1, 2\}$, the curve $C_i$ has an ordinary triple point $Q_i$, with $Q_1 = (0:0:1)$, $Q_2 = (0:1:1)$. The curve $C_3$ has an ordinary double point at $Q_3 = (1:1:1)$. For all $i$, the point $Q_i$ lies on the branch locus $B_i = B_{X_i}$.

We will see later that the curve $C_i$ intersects the branch locus $B_i$ with even multiplicity everywhere. Of course, one could check this directly as well using the polynomial $h_i$. In fact, had we *defined* $C_i$ by the vanishing of $h_i$, then one would easily check that $C_i$ satisfies the conditions of part (2) of Corollary 1.3, which gives an alternative proof of unirationality of $X_i$ without the need of the explicit morphism $\rho_i$ (see Example 5.16). Indeed, in practice we first found the curves $C_1$, $C_2$, and $C_3$, and then constructed the parametrisations $\rho_1, \rho_2, \rho_3$, which allow for the more direct proof that we gave of Theorem 3.3.

The following section gives an alternative description of Yu. Manin's original construction of curves on del Pezzo surfaces of degree 2 and the generalisation of C. Salgado, D. Testa, and A. Várilly-Alvarado. The benefit of this description is that it allows us to write down the curves explicitly, in full generality. This will be used in Section 6. It also shows that quartic plane curves with a triple point and cubic points with a double point play an essential role in Manin's construction as well.

## 4. Manin's construction

Let $k$ be a field and $X$ a del Pezzo surface of degree 2 over $k$. Let $\pi \colon X \to \mathbb{P}^2$ be the associated double covering map, with ramification locus $R \subset X$, and associated involution $\iota \colon X \to X$. Let $P \in X(k)$ be a point. For this entire section, we will assume that $P$ does not lie on $R$. Yu. Manin proves unirationality of $X$ under the extra condition that $P$ does not lie on any of the 56 exceptional curves of $X$ (see [Man86, Theorem 29.4]). His proof relies on the following construction of a non-constant morphism $\rho \colon \mathbb{P}^1 \to X$. Let $\mathrm{bl}_P \colon \tilde{X} \to X$ be the blow-up of $X$ at $P$. Assume for the remainder of this paragraph that $P$ does not lie on an exceptional curve of $X$. Then $\tilde{X}$ is a del Pezzo surface of degree 1 (see [STVA14, Corollary 14]), on which there are 240 exceptional curves, splitting up naturally in 120 pairs $(E, E')$ with $E + E'$ linearly equivalent to $-2K_{\tilde{X}}$, where $K_{\tilde{X}}$ is a canonical divisor on $\tilde{X}$. Let $E'_0 \subset \tilde{X}$ be the exceptional curve on $\tilde{X}$ that is paired with the exceptional curve $E_0$ above $P$. Then any isomorphism $\mathbb{P}^1 \to E'_0$ composed with the restriction of $\mathrm{bl}_P$ to $E'_0$ yields a morphism $\rho \colon \mathbb{P}^1 \to X$. Unirationality of $X$ then follows from applying the same construction to the generic point of $\mathrm{bl}_P(E'_0)$ instead of $P$ (see [Man86, Theorem 29.4]).

C. Salgado, D. Testa, and A. Várilly-Alvarado generalise Manin's construction to the case that $P$ is not a generalised Eckardt point, that is, $P$ lies on at most three exceptional curves. They show that in this case there is still a unique effective divisor on $\tilde{X}$ that is linearly equivalent to $-2K_{\tilde{X}} - E_0 \sim -2\,\mathrm{bl}_P^* K_X - 3E_0$, and this divisor contains a component that is birational to $\mathbb{P}^1$ (see [STVA14, Theorem 16 and its proof]).

We give a different description of this rational component by generalising a well-known alternative description of Manin's construction in terms of an elliptic fibration associated to the anti-canonical map $\tilde{X} \dashrightarrow \mathbb{P}^1$. We drop all assumptions on the number of exceptional curves going through $P$. Let $\mathrm{bl}_{P,\iota(P)} \colon \mathrm{Bl}_{P,\iota(P)} X \to X$ denote the blow-up of $X$ at the points $P$ and $\iota(P)$, which factors as the composition of $\mathrm{bl}_P$ and the blow-up $\mathrm{bl}_{\iota(P)}$ of $\tilde{X}$ at the point corresponding to $\iota(P)$. Let $E_P$ and $E_{\iota(P)}$ be the exceptional curves above $P$ and $\iota(P)$, respectively. Let $\mathrm{bl}_Q \colon \mathrm{Bl}_Q \mathbb{P}^2 \to \mathbb{P}^2$ denote the blow-up of $\mathbb{P}^2$ at $Q$. The map $\pi \colon X \to \mathbb{P}^2$ induces a morphism $\pi' \colon \mathrm{Bl}_{P,\iota(P)} X \to \mathrm{Bl}_Q \mathbb{P}^2$. Let $\mathfrak{L}_Q$ denote the line in the dual of $\mathbb{P}^2$ consisting of all lines $L \subset \mathbb{P}^2$ going through $Q$, and note that $\mathfrak{L}_Q$ is isomorphic to $\mathbb{P}^1$. We can identify the blow-up $\mathrm{Bl}_Q \mathbb{P}^2$ with the subvariety of $\mathbb{P}^2 \times \mathfrak{L}_Q$ consisting of pairs $(R, L)$ with $R \in L$. Let $\psi \colon \mathrm{Bl}_{P,\iota(P)} X \to \mathfrak{L}_Q$ be the composition of $\pi'$ with the projection $\mathrm{Bl}_Q \mathbb{P}^2 \to \mathfrak{L}_Q$. The fiber of $\psi$ above a line $L$ through $Q$ is isomorphic to the pull back $\pi^*(L)$ on $X$, which is linearly equivalent to $-K_X$. By the adjunction formula we have

$$2g(\pi^*(L)) - 2 = \pi^*(L) \cdot (\pi^*(L) + K_X) = 0,$$

from which it follows that $g(\pi^*(L)) = 1$, so $\psi$ is an elliptic fibration. The morphism $\pi'$ maps the exceptional curves $E_P$ and $E_{\iota(P)}$ isomorphically to the exceptional curve on $\mathrm{Bl}_Q \mathbb{P}^2$ above $Q$, which in turn projects isomorphically to $\mathfrak{L}_Q$. Hence, we obtain two sections $\sigma_P$ and $\sigma_{\iota(P)}$ of $\psi$, mapping to $E_P$ and $E_{\iota(P)}$, respectively. The rational map $\psi \circ \mathrm{bl}_{\iota(P)}^{-1} \colon \tilde{X} \dashrightarrow \mathfrak{L}_Q \cong \mathbb{P}^1$ is associated to the linear system $|-K_{\tilde{X}}|$ on $\tilde{X}$, which has a unique base point, namely the point associated to $\iota(P)$. It also factors as the composition of $\pi \circ \mathrm{bl}_P$ and the projection $p_Q \colon \mathbb{P}^2 \dashrightarrow \mathfrak{L}_Q$ away from $Q$. These maps are summarised in the following commutative diagram.

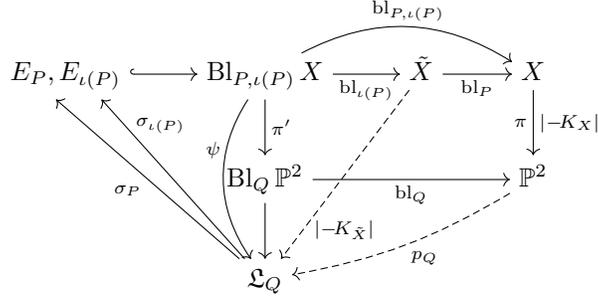

We fix the section $\sigma_{\iota(P)}$ of $\psi$ as the neutral element of the Mordell–Weil group $(\mathrm{Bl}_{P,\iota(P)} X)(\mathfrak{L}_Q)$ of sections of $\psi$ and we set

(2) $$\rho = \mathrm{bl}_{P,\iota(P)} \circ (-\sigma_P) \colon \mathfrak{L}_Q \to X.$$

The following Proposition is well known among experts and implies that there is no conflict with the definition of $\rho$ in Manin's construction above, at least after an appropriate identification of $\mathfrak{L}_Q$ with $\mathbb{P}^1$.

**Proposition 4.1.** *If $P$ does not lie on an exceptional curve, then $\mathrm{bl}_{\iota(P)} \circ (-\sigma_P) \colon \mathfrak{L}_Q \to \tilde{X}$ induces an isomorphism from $\mathfrak{L}_Q$ to the exceptional curve $E'_0 \subset \tilde{X}$ of Manin's construction.*

*Proof.* Given that $P$ does not lie on the ramification curve $R$, the assumption that $P$ does not lie on an exceptional curve implies that $\tilde{X}$ is a del Pezzo surface of degree 1 (see [STVA14, Corollary 13]). Therefore, the lemma follows from the following general statement. For any del Pezzo surface $Y$ of degree 1, the anti-canonical map $Y \dashrightarrow \mathbb{P}^1$ induces an elliptic fibration of the blow-up $\tilde{Y}$ of $Y$ in the unique base point of the linear system $|-K_Y|$, with the exceptional curve as the zero section. Multiplication by $-1$ induces an involution of $Y$ that switches any two exceptional curves $E$ and $E'$ for which $E + E'$ is linearly equivalent to $-2K_Y$. □

We continue to investigate the general case, in which $P$ may lie on any number of exceptional curves.

**Proposition 4.2.** *The following statements are equivalent.*
  (1) *The map $\rho \colon \mathfrak{L}_Q \to X$ is not constant.*
  (2) *The section $\sigma_P$ does not have order $2$.*
  (3) *The composition $p_Q \circ \pi \circ \rho \colon \mathfrak{L}_Q \to X \to \mathbb{P}^2 \dashrightarrow \mathfrak{L}_Q$ is a well-defined rational map that extends to the identity on $\mathfrak{L}_Q$.*

*Proof.* The map $\rho$ is constant if and only if the image of $-\sigma_P$ is contained in $E_P$ or $E_{\iota(P)}$. Since $\sigma_P, \sigma_{\iota(P)}$, and $-\sigma_P$ are all sections of $\varphi$, this happens if and only if $-\sigma_P = \sigma_P$ or $-\sigma_P = \sigma_{\iota(P)}$. The latter case does not happen as $\sigma_{\iota(P)} = 0$ in the group $(\mathrm{Bl}_{P,\iota(P)} X)(\mathfrak{L}_Q)$. Hence, the map $\rho$ is constant if and only if $\sigma_P = -\sigma_P$, which proves the equivalence of (1) and (2). Moreover, if $\rho$ is constant, then it sends $\mathfrak{L}_Q$ to $P$ and the image of $\pi \circ \rho$ is then $Q$. Since $\pi$ is finite, the map $\rho$ is non-constant if and only if $\pi \circ \rho$ is non-constant. As the base locus of $p_Q$ consists of just the one point $Q$, the map $\pi \circ \rho$ is non-constant if and and only if the composition $p_Q \circ \pi \circ \rho$ is a well-defined rational map. If this is the case, then by the commutivity of the diagram above, this composition coincides with $\psi \circ (-\sigma_P) = \mathrm{id}_{\mathfrak{L}_Q}$. This proves that (3) is equivalent as well. □

We will see later, in Theorem 4.13, that the conditions in Proposition 4.2 are satisfied if and only if $P$ is not a generalised Eckardt point.

We set $B = \pi(R)$, so $B$ is the branch curve of $\pi$. Recall that the fiber of $\psi$ above a line $L \in \mathfrak{L}_Q$ is isomorphic to the pull back $\pi^*(L)$, which is singular if and only if $L$ is tangent to $B$, and reducible if and only if $L$ is a bitangent to $B$. In the latter case, the fiber has Kodaira type $I_2$, with the two components being the strict transforms of the two exceptional curves of $X$ above $L$. We now state the first alternative description of $\rho$.

**Proposition 4.3.** *Let $L \in \mathfrak{L}_Q$ be a line that is not tangent to the branch curve $B$. Then $\rho(L)$ is the unique point on $\pi^{-1}(L)$ that is linearly equivalent with the divisor $2(\iota(P)) - (P)$ on $\pi^{-1}(L)$.*

*Proof.* Since $L$ is not tangent to $B$, the curve $C = \pi^{-1}(L)$ is smooth and has genus 1. By the definition of $\rho$, we find $\rho(L) = -P$ in the group $C(k)$ with $\iota(P)$ as the neutral element. The statement therefore follows from a standard result on elliptic curves (see [Sil09, Corollary III.3.5]). □

The following lemma will allow us to describe $\rho(L)$ more explicitly, also in the more general case that $L$ is not a *bi*-tangent to $B$.

**Lemma 4.4.** *Let $x, y, w$ denote the coordinates of the weighted projective space $\mathbb{P}(1,1,2)$ over $k$. Set $\Omega = (0:0:1)$ and $U = \mathbb{P}(1,1,2) - \{\Omega\}$, and let $\alpha\colon U \to \mathbb{P}^1$ be the projection away from $\Omega$. Let $C \subset \mathbb{P}(1,1,2)$ be a reduced, not-necessarily-irreducible curve given by $w^2 + fw = g$ for some homogeneous polynomials $f, g \in k[x, y]$ of degree 2 and 4, respectively. Let $S \in C(k)$ be a rational point where the restriction $\alpha|_C\colon C \to \mathbb{P}^1$ does not ramify. Then there is a unique section $\tau\colon \mathbb{P}^1 \to U$ of $\alpha$ of which the image intersects $C$ at $S$ with multiplicity at least 3. The following hold.*
  (1) *The image $\tau(\mathbb{P}^1)$ has degree 2 in $\mathbb{P}(1,1,2)$.*
  (2) *If $C$ is reducible, then $\tau$ induces an isomorphism from $\mathbb{P}^1$ to the irreducible component of $C$ containing $S$.*
  (3) *If $C$ is irreducible, then $\tau(\mathbb{P}^1)$ intersects $C$ in a unique fourth point. If this point is different from $S$, then the corresponding intersection multiplicity is 1.*

*Proof.* Let $\tau\colon \mathbb{P}^1 \to U$ be a section of $\alpha$. Then there are homogeneous polynomials $q, r \in k[x, y]$ with no factor of $q$ dividing $r$ doubly, such that $\tau$ sends $(x:y)$ to $(xq(x,y) : yq(x,y) : r(x,y))$. Any point of $\mathbb{P}^1$ on which $q$ vanishes maps to $\Omega \notin U$, so we conclude that $q$ is a constant and thus $r$ has degree 2. Without loss of generality we assume $q = 1$. The image $\tau(\mathbb{P}^1)$ is given by $w = r(x, y)$, so $\tau(\mathbb{P}^1)$ has degree 2 as well, which proves (1). After switching $x$ and $y$ if necessary, we may assume that $S$ lies in the affine part given by $y \neq 0$. Set $v = w/y^2$ and $t = x/y$ and let $a, b, c \in k$ be such that $r(x, y) = ax^2 + bxy + cy^2$. Then on the curve $\tau(\mathbb{P}^1)$, given by $v = r(t, 1)$, we have

$$v(S) = r(t,1)(S), \qquad \frac{dv}{dt}(S) = \frac{dr(t,1)}{dt}(S), \qquad \text{and} \qquad \frac{d^2v}{dt^2}(S) = \frac{d^2r(t,1)}{dt^2}(S).$$

Hence, the condition that $\tau(\mathbb{P}^1)$ intersects $C$ with multiplicity at least 3 at $S$ is equivalent with these equalities holding on $C$ as well. This gives three linearly independent linear equations in $a$, $b$, and $c$, so there are unique $a, b, c \in k$ for which $\tau(\mathbb{P}^1)$ intersects $C$ at $S$ with multiplicity at least 3. Suppose that $\tau$ is this section. Given that $\alpha|_C$ does not ramify at $S$, there is a unique irreducible component $C_0$ of $C$ containing $S$. If $C_0$ and $\tau(\mathbb{P}^1)$ do not coincide, then the total number of their intersection points, counted with multiplicities, is $(\deg C_0) \cdot (\deg \tau(\mathbb{P}^1))/(1 \cdot 1 \cdot 2) = \deg C_0$. Since $C$ is a double cover of $\mathbb{P}^1$, either $C$ is irreducible, in which case $C_0 = C$ has degree 4, so there is a unique fourth intersection point (proving (3)), or $C$ has two isomorphic irreducible components of degree 2. In the latter case, the intersection multiplicity at $S$ exceeds the total intersection number $C_0 \cdot \tau(\mathbb{P}^1)$. From this contradiction, we conclude that $C_0$ and $\tau(\mathbb{P}^1)$ do coincide; as $\tau$ is a section of $\alpha|_{C_0}$, we find that $\tau$ induces an isomorphism from $\mathbb{P}^1$ to $C_0$, which proves (2). □

*Remark* 4.5. In the case that $C$ in Lemma 4.4 is smooth, the Riemann-Roch Theorem can be used to give a more efficient proof.

*Remark* 4.6. Assume the same setting as in Lemma 4.4. We will make the section $\tau$ explicit. Write $f$ and $g$ as

$$f = f_2 x^2 + f_1 xy + f_0 y^2 \qquad \text{and} \qquad g = g_4 x^4 + g_3 x^3 y + g_2 x^2 y^2 + g_1 xy^3 + g_0 y^4.$$

Without loss of generality, we assume $S = (0:1:w_0)$ for some $w_0 \in k$ with $w_0^2 + f_0 w_0 = g_0$. Write $\delta = f_0 + 2w_0$. The fact that $\alpha|_C$ does not ramify at $S$ is equivalent with $\delta \neq 0$. As in the proof of Lemma 4.4, we set $v = w/y^2$ and $t = x/y$, and we let $a, b, c \in k$ be such that $\tau$ sends $(x:y)$ to $(x:y:r(x,y))$ with $r(x,y) = ax^2 + bxy + cy^2$. From $\tau((0:1)) = S$, we find $c = w_0$. In terms of the affine coordinates $t, v$, the intersection points of $\tau(\mathbb{P}^1)$ and $C$ correspond to the roots of $F(t) = r(t,1)^2 + f(t,1)r(t,1) - g(t,1)$. Write

$$F(t) = F_4 t^4 + F_3 t^3 + F_2 t^2 + F_1 t + F_0,$$

with $F_i \in k$. Then we have

$$\begin{aligned} F_0 &= 0, \\ F_1 &= \delta b + w_0 f_1 - g_1, \\ F_2 &= \delta a + b^2 + f_1 b + w_0 f_2 - g_2, \\ F_3 &= 2ab + f_1 a + f_2 b - g_3, \\ F_4 &= a^2 + f_2 a - g_4. \end{aligned}$$

The intersection multiplicity of $\tau(\mathbb{P}^1)$ and $C$ at $S$ equals the multiplicity of $t = 0$ as a root of $F(t)$. Since this multiplicity is at least 3, we find $F_1 = F_2 = 0$. This determines $a$ and $b$ uniquely, as we have $\delta \neq 0$. If, for these $a, b$, we have $F_3 = F_4 = 0$ as well, then $\tau(\mathbb{P}^1)$ is contained in $C$, so $C$ consists of two components. Otherwise, the fourth intersection point corresponds to $t = -F_3/F_4$, so it is $(-F_3 : F_4 : aF_3^2 - bF_3F_4 + w_0F_4^2)$.

We identify $X$ with a model in $\mathbb{P}$ given by $w^2 + fw = g$ for some homogeneous polynomials $f, g \in k[x, y, z]$ of degrees 2 and 4, respectively. Set $\Omega = (0 : 0 : 0 : 1) \in \mathbb{P}$ and $U = \mathbb{P} - \{\Omega\}$. The projection $\varphi \colon U \to \mathbb{P}^2$ onto the first three coordinates, that is, the projection away from $\Omega$, restricts to the map $\pi \colon X \to \mathbb{P}^2$.

**Proposition 4.7.** *Let $L \in \mathfrak{L}_Q$ be a line through $Q$. Then there is a unique morphism $\tau \colon L \to U$ with $\varphi \circ \tau = \mathrm{id}_L$ of which the image intersects $X$ at $P$ with multiplicity at least 3. The following hold.*

(1) *The image $\tau(L)$ has degree 2 in $\mathbb{P}$.*
(2) *If $L$ is bitangent to the branch curve $B$, then $\tau$ induces an isomorphism from $L$ to the unique exceptional curve above $L$ that contains $P$.*
(3) *If $L$ is not bitangent to the branch curve $B$, then $\tau(L)$ intersects $X$ in a unique fourth point, which is $\rho(L)$. If $\rho(L) \neq P$, then the corresponding intersection multiplicity is 1.*

*Proof.* The inverse image $X_L = \pi^{-1}(L)$ is contained in the inverse image $U_L = \varphi^{-1}(L)$. Note that every morphism $\tau \colon L \to U$ with $\varphi \circ \tau = \mathrm{id}_L$ factors through $U_L$, and the intersection multiplicity of $\tau(\mathbb{P})$ with $X$ at $P$ equals the intersection multiplicity of $\tau(\mathbb{P})$ with $X_L$ at $P$. Since $U_L \cup \{\Omega\}$ is isomorphic to the weighted projective space $\mathbb{P}(1, 1, 2)$, the unique existence of $\tau$ follows from Lemma 4.4, applied to $C = X_L$.

$$\begin{array}{ccc} X_L & \hookrightarrow U_L \xrightarrow{\tau} & L \\ \downarrow & \downarrow & \downarrow \\ X & \hookrightarrow U \xrightarrow{\varphi} & \mathbb{P}^2 \\ & \underset{\pi}{\longrightarrow} & \end{array}$$

The same lemma also implies (1) and (2). Suppose $L$ is not bitangent to $B$. Then by Lemma 4.4, the image $\tau(L)$ intersects $X_L$, and thus $X$, in a unique fourth point, say $T$, and the corresponding intersection multiplicity is 1 if $T \neq P$. By (1) of Lemma 4.4, the divisor $(T) + 3(P)$ is linearly equivalent to the intersection of $X_L$ with a curve of degree 2 in $U_L \cup \{\Omega\} \cong \mathbb{P}(1, 1, 2)$. Since $(P) + (\iota(P))$ is the intersection of $X_L$ with a curve of degree 1, and all curves in $\mathbb{P}(1, 1, 2)$ of the same degree are linearly equivalent, we find that the divisors $(T) + 3(P)$ and $2(P) + 2(\iota(P))$ on $X_L$ are linearly equivalent. From Proposition 4.3 we find $T = \rho(L)$, which finishes the proof of (3). □

From now on, for any $L \in \mathfrak{L}_Q$, we let

$$\tau_L \colon L \to U$$

denote the map associated to $L$ as in Proposition 4.7. Let the map $\nu\colon \mathrm{Bl}_Q \mathbb{P}^2 \to U$ be given by sending the pair $(S, L) \in \mathrm{Bl}_Q \mathbb{P}^2 \subset \mathbb{P}^2 \times \mathfrak{L}_Q$ with $S \in L$ to $\tau_L(S)$.

**Proposition 4.8.** *The map $\nu\colon \mathrm{Bl}_Q \mathbb{P}^2 \to U$ is a morphism.*

*Proof.* Since $\nu$ is well defined globally, it suffices to show that $\nu$ is a morphism locally. Let $T \in \mathbb{P}^2$ be a point not equal to $Q$, and let $V_T \subset \mathrm{Bl}_Q \mathbb{P}^2$ be the open subset of all $(S, L)$ with $T \notin L$. Let $\mathfrak{L}_T$ denote the line in the dual of $\mathbb{P}^2$ consisting of all lines $L \subset \mathbb{P}^2$ going through $T$. Let $\mathfrak{L}_Q^\circ \subset \mathfrak{L}_Q$ denote the subset of all $L \in \mathfrak{L}_Q$ with $T \notin L$. For each $L \in \mathfrak{L}_Q^\circ$, there is an isomorphism $\mu_L \colon L \to \mathfrak{L}_T$ that sends $S \in L$ to the line through $T$ and $S$. The map $\nu_0 \colon V_T \to \mathfrak{L}_T \times \mathfrak{L}_Q^\circ$ sending $(S, L)$ to $(\mu_L(S), L)$ is an isomorphism. Let $\mathrm{Mor}(\mathfrak{L}_T, U)$ denote the scheme of morphisms from $\mathfrak{L}_T$ to $U$. Define $\alpha\colon \mathfrak{L}_Q^\circ \to \mathrm{Mor}(\mathfrak{L}_T, U)$ to be the map that sends $L$ to $\tau_L \circ (\mu_L)^{-1}$. We may identify $\mathfrak{L}_T$ with $\mathbb{P}^1$, and write $\mathcal{O}_{\mathfrak{L}_T}(n)$ for the sheaf on $\mathfrak{L}_T$ corresponding to $\mathcal{O}_{\mathbb{P}^1}(n)$, for any integer $n$. For any $L \in \mathfrak{L}_Q^\circ$, it follows from Remark 4.6 that there exist elements $\xi, \eta, \zeta \in \Gamma(\mathcal{O}_{\mathfrak{L}_T}(1))$ and $\omega \in \Gamma(\mathcal{O}_{\mathfrak{L}_T}(2))$, such that $\alpha(L)\colon \mathfrak{L}_T \to U \subset \mathbb{P}$ is given by $S \mapsto (\xi(S) : \eta(S) : \zeta(S) : \omega(S))$. Since we can do this for the generic point of $\mathfrak{L}_Q^\circ$ as well, we find that $\alpha$ is a rational map, and as $\mathfrak{L}_Q^\circ$ is a smooth curve, the map $\alpha$ is in fact a morphism. On $V_T$, the map $\nu$ factors as

$$V_T \xrightarrow{\nu_0} \mathfrak{L}_T \times \mathfrak{L}_Q^\circ \xrightarrow{(\mathrm{id}_{\mathfrak{L}_T}, \alpha)} \mathfrak{L}_T \times \mathrm{Mor}(\mathfrak{L}_T, U) \xrightarrow{\mathrm{ev}} U$$

where the last morphism is the evaluation map given by $\mathrm{ev}(L, \beta) = \beta(L)$. Hence, the restriction of $\nu$ to $V_T$ is a morphism. Since $T$ was arbitrary, we conclude that $\nu$ is a morphism as well. $\square$

The following proposition shows that we can combine all the maps $\tau_L \colon L \to U$ to one map $\mathbb{P}^2 \to U$.

**Proposition 4.9.** *The map $\nu \colon \mathrm{Bl}_Q \mathbb{P}^2 \to U$ factors through a unique morphism $\xi \colon \mathbb{P}^2 \to U$. Moreover, the following hold.*

(1) *For any $L \in \mathfrak{L}_Q$, the restriction $\xi|_L$ of $\xi$ to $L$ equals $\tau_L$.*
(2) *The map $\xi$ is a section of the map $\varphi \colon U \to \mathbb{P}^2$ with $\xi(Q) = P$.*
(3) *The image $\xi(\mathbb{P}^2)$ has degree 2.*
(4) *The intersection $X \cap \xi(\mathbb{P}^2)$ is the union of $\rho(\mathfrak{L}_Q)$ and the exceptional curves on $X$ that contain $P$.*

*Proof.* By Zariski's Main Theorem (see [Har77, Corollary III.11.4 and its proof]), we have an isomorphism $(\mathrm{bl}_Q)_* \mathcal{O}_{\mathrm{Bl}_Q \mathbb{P}^2} \to \mathcal{O}_{\mathbb{P}^2}$. Note also that the map $\nu$ contracts the exceptional curve above $Q$ to the point $P$. Therefore, it follows from [Gro61, Lemme 8.11.1] that $\nu$ factors as the composition of $\mathrm{bl}_Q$ and a unique morphism $\xi \colon \mathbb{P}^2 \to U$. Let $L$ be a line through $Q$ and $S$ a point on $L$. Then we have $\xi(S) = \nu(S, L) = \tau_L(S)$, which proves (1). It also follows that we have $\varphi(\xi(S)) = \varphi(\tau_L(S)) = S$, so $\xi$ is a section of $\varphi$. From $\tau_L(Q) = P$, we find $\xi(Q) = P$, which proves (2). The line $\varphi^{-1}(S)$ intersects $Y = \xi(\mathbb{P}^2)$ in the same number of points as it intersects $\tau_L(L)$, so the degree of $Y$ equals 2 by Proposition 4.7. This proves (3). It follows from the same proposition that the intersection of $\tau_L(L) = \xi(L)$ with $X$ is equal to $\{P, \rho(L)\}$ if $L$ is not bitangent to $B$, and equal to the exceptional curve above $L$ that contains $P$ if $L$ is bitangent to $B$. As $Y = \xi(\mathbb{P}^2)$ is the union of the images $\xi(L)$ as $L$ runs through $\mathfrak{L}_Q$, the last statement follows as well. $\square$

*Remark* 4.10. We can make the morphism $\xi \colon \mathbb{P}^2 \to U$ completely explicit. Let $X \subset \mathbb{P}$ be given by $w^2 + fw = g$ with $f, g \in k[x, y, z]$ homogeneous of degree 2 and 4, respectively. Write $f$ and $g$ as

$$f = f_2 + f_1 z + f_0 z^2 \quad \text{and} \quad g = g_4 + g_3 z + g_2 z^2 + g_1 z^3 + g_0 z^4,$$

with $f_i, g_i \in k[x, y]$ homogeneous of degree $i$. Without loss of generality, we assume $Q = (0 : 0 : 1)$, that is, $P = (0 : 0 : 1 : w_0)$ for some $w_0 \in k$ with $w_0^2 + f_0 w_0 = g_0$. Write $\delta = 2w_0 + f_0$. The fact that $\pi$ does not ramify at $P$ is equivalent with $\delta \neq 0$. Let $\xi' \colon \mathbb{P}^2 \to U$ be the morphism that sends $(x : y : z)$ to $(x : y : z : w)$ with $w$ satisfying

(3) $$\delta^3 w = \delta^2(g_2 + g_1 z + 2g_0 z^2 - w_0 f) - g_1^2 + f_1(g_0 f_1 - f_0 g_1).$$

We claim $\xi = \xi'$. To prove the claim, one first checks that this $w$ equals $a + bz + w_0 z^2$ with $a$ and $b$ satisfying

(4) $$\delta b + w_0 f_1 - g_1 = 0 \quad \text{and} \quad \delta a + b^2 + f_1 b + w_0 f_2 - g_2 = 0.$$

Now take any $\lambda \in k$. By projecting away from the $y$-coordinate, we identify the line $L \subset \mathbb{P}^2$ and the hyperplane $H \subset \mathbb{P}$, both given by $y = \lambda x$, with $\mathbb{P}^1$ and $\mathbb{P}(1,1,2)$, respectively. The intersection $H \cap X$ gets identified with a curve $C \in \mathbb{P}(1,1,2)$. Under these identifications, the equations (4) correspond to the vanishing of $F_1$ and $F_2$ as in Remark 4.6. This implies that the restriction $\xi'|_L$ of $\xi'$ to $L$ corresponds to the map $\tau_L \colon L \to U$. From part (1) of Proposition 4.9 we therefore conclude $\xi' = \xi$. This explicit description of $\xi$ could have also been used to give a less-enlightening proof of Propositions 4.8 and 4.9.

*Remark* 4.11. Assume the characteristic of $k$ is not 2. We have assumed that $P$ does not lie on the ramification curve of $\pi$, but we may deform the surface $X$ and the point $P$ such that in the limit, the point does lie on the ramification curve. In the limit we have $\delta = 0$ and $f_0 = -2w_0$ and $g_0 = w_0^2 + f_0 w_0 = -w_0^2$, which implies that the limit of $\xi(\mathbb{P}^2)$, which is given by (3), is the non-reduced surface in $\mathbb{P}$ given by the vanishing of

$$-g_1^2 + f_1(g_0 f_1 - f_0 g_1) = -(g_1 - w_0 f_1)^2.$$

The intersection of the surface given by $g_1 - w_0 f_1 = 0$ with $U = \mathbb{P} - \{\Omega\}$ is exactly the pull back under the projection $\varphi \colon U \to \mathbb{P}^2$ of the tangent line $L$ to $B$ at $Q$, so we find that, as a divisor, the limit of the intersection $X \cap \xi(\mathbb{P}^2)$ is $2\pi^*(L)$.

It follows from Propositions 4.12 and 4.14 below that the divisor $\mathrm{bl}_P^*(\xi(\mathbb{P}^2) \cap X) - 3E_P$ on $\tilde{X}$ corresponds with the effective divisor denoted by $D$ in the proof of [STVA14, Theorem 16]. We will mimic part of this proof, working on $X$ instead of $\tilde{X}$, to show the divisor $\xi(\mathbb{P}^2) \cap X$ is reduced.

**Proposition 4.12.** *Set $Y = \xi(\mathbb{P}^2) \in \mathrm{Div}\,\mathbb{P}$ and let $D \in \mathrm{Div}\,X$ be the pull back of $Y$ to $X$. Then $D$ is linearly equivalent with $-2K_X$. We have $D = D_1 + \sum_E E$, where the sum runs over all exceptional curves $E$ containing $P$, and where $D_1 = \rho(\mathfrak{L}_Q)$ if $\rho$ is not constant and $D_1 = 0$ otherwise.*

*Proof.* By Proposition 4.9, the image $Y$ of $\xi$ is a hypersurface in $\mathbb{P}$ of degree 2, which implies the first statement. Let $\mathcal{E}$ denote the set of all exceptional curves $E$ containing $P$. By Proposition 4.9, we may write $D = D_1 + \sum_E n_E E$, where the sum runs over all $E \in \mathcal{E}$, and where $D_1$ is effective and supported on the integral curve $\rho(\mathfrak{L}_Q)$ if $\rho$ is not constant and $D_1 = 0$ otherwise. If $\rho$ is not constant, then $D_1 = m \cdot \rho(\mathfrak{L}_Q)$ for some positive integer $m$. If we had $m > 1$, that is, the intersection $Y \cap X$ is nonreduced along $\rho(\mathfrak{L}_Q)$, then for any line $L \subset \mathbb{P}^2$ through $Q$, we find, after intersecting with the hyperplane $H_L = \overline{\varphi^{-1}(L)} \subset \mathbb{P}$, that the intersection of $\tau_L(L) = Y \cap H_L$ with $X$ is not smooth at the point $\rho(L)$. This contradicts Proposition 4.7, part (3), when $\rho(L) \neq P$, so we conclude $m = 1$ and thus $D_1 = \rho(\mathfrak{L}_Q)$ in the case that $\rho$ is not constant.

It remains to prove $n_E = 1$ for all $E \in \mathcal{E}$. From Proposition 4.9 we find $n_E \geq 1$ for all $E \in \mathcal{E}$. Suppose there is an $E_0 \in \mathcal{E}$ with $n_{E_0} > 1$. Then $D - 2E_0$ is an effective divisor that is linearly equivalent with $-2K_X - 2E_0$, and thus with $2\iota(E_0)$. Since these divisors have negative self-intersection, their linear system contains only one effective divisor, so we find $D - 2E_0 = 2\iota(E_0)$. As $P$ is not contained in the ramification locus, the exceptional curve $\iota(E_0)$ does not contain $P$, so it is not contained in $\mathcal{E}$. Hence, we find $D_1 = 2\iota(E_0)$, which implies that $\rho$ is not constant and $\rho(\mathfrak{L}_Q) = \iota(E_0)$. But then the equality $D_1 = 2\iota(E_0)$ contradicts $m = 1$, so we conclude that we have $n_E = 1$ for all $E \in \mathcal{E}$. $\square$

We are now ready to prove that the conditions of Proposition 4.2 hold if and only if $P$ is not a generalised Eckardt point (cf. [STVA14, Theorem 16]).

**Theorem 4.13.** *The morphism $\rho \colon \mathfrak{L}_Q \to X$ is constant if and only if $P$ is a generalised Eckardt point.*

*Proof.* Let $\mathcal{E}$ denote the set of all exceptional curves containing $P$, and set $n = \#\mathcal{E}$. Let $D$ and $D_1$ be as in Proposition 4.12, so that $D_1$ is linearly equivalent with $-2K_X - \sum_E E$, where the sum runs over all $E \in \mathcal{E}$. Then we have $-K_X \cdot D_1 = 4 - n$. As $-K_X$ is ample and $D_1$ is effective, we have $-K_X \cdot D_1 \geq 0$ with equality if and only if $D_1 = 0$. We conclude that $P$ is a generalised Eckardt point, that is, we have $n = 4$, if and only if $D_1 = 0$, that is, $\rho$ is constant. $\square$

The following proposition describes the image of $\rho$ more precisely in the case that $P$ is not a generalised Eckardt point.

**Proposition 4.14.** *Suppose $P$ is not a generalised Eckardt point. Let $n$ be the number of lines in $\mathfrak{L}_Q$ that are bitangent to the branch curve $B$. Then $\pi$ induces an isomorphism from $\rho(\mathfrak{L}_Q)$ to its image $C$. Furthermore, the curve $C$ has degree $4 - n$ and $Q$ is a point of multiplicity $3 - n$ on $C$.*

*Proof.* By Theorem 4.13, the morphism $\rho$ is not constant, so $C$ is an irreducible curve. Since $P$ is not a ramification point of $\pi$, the number of exceptional curves through $P$ is $n$. Let $D$ and $D_1 = \rho(\mathfrak{L}_Q)$ be as in Proposition 4.12. Note that the same proposition tells us that $D$ and $D_1$ are reduced. Since $\xi$ is a section of $\varphi \colon U \to \mathbb{P}^2$, and $\varphi$ restricts to $\pi \colon X \to \mathbb{P}^2$, the map $\pi$ induces an isomorphism from $D = \xi(\mathbb{P}^2) \cap X$ onto its image. In particular, it induces an isomorphism $D_1 \to C$, so we find $\pi_*(D_1) = C$. Let $M \subset \mathbb{P}^2$ be any line not containing $Q$. Then the projection formula (see [Kol96, Proposition VI.2.11]) implies
$$\deg C = C \cdot M = \pi_*(D_1) \cdot M = D_1 \cdot \pi^*(M) = D_1 \cdot (-K_X) = 4 - n.$$
Let $m$ be the multiplicity of $Q$ on $C$. By Proposition 4.2, the projection map $p_Q \colon \mathbb{P}^2 \dashrightarrow \mathfrak{L}_Q$ induces a birational map $C \dashrightarrow \mathfrak{L}_Q$. This implies $m = \deg C - 1 = 3 - n$. $\square$

*Remark* 4.15. Proposition 4.14 shows that if $P$ is not a generalised Eckardt point, then $C = \pi(\rho(\mathfrak{L}_Q))$ is a curve in $\mathbb{P}^2$ with a point of multiplicity $\deg C - 1$ off the branch locus $B$. Moreover, in the next section, we will see that $C$ intersects $B$ with even intersection multiplicity everywhere. We will also see that outside characteristic 2, all such curves pull back to rational curves on $X$ and imply unirationality of $X$. See Example 5.14.

## 5. Proof of the second main theorem

From now on, we assume that the characteristic of $k$ is not 2. We continue the notation of the previous section, but we drop the assumption that $P$ does not lie on the ramification locus $R$. In the previous section, we described the image in $\mathbb{P}^2$ of the rational curve on $X$ that is obtained through Manin's construction. In this section we will conversely describe all curves $C \subset \mathbb{P}^2$ of which the pull back to $X$ contains a component of genus 0 (cf. Remark 5.8). Whenever such a component is rational over $k$, the surface $X$ is unirational by Lemma 3.2. The first subsection gives a description in terms of the normalisation of $C$, while the second subsection gives a description in terms of $C$ itself, under a mild hypothesis.

### 5.1. The pull back of a plane curve in terms of its normalisation.
Let $\overline{k}$ denote an algebraic closure of $k$. Given a reduced curve $C$ over $k$, the normalisation map $\vartheta \colon \tilde{C} \to C$ is unique up to isomorphism; the curve $\tilde{C}$ is regular and both $\vartheta$ and $\tilde{C}$ are called the normalisation of $C$. For more details, see [Mum99, Theorem III.8.3] for the case that $C$ is irreducible. For the general case, take the disjoint union of the normalisations of the irreducible components. If $X$ is a variety over $k$ and $\ell$ is a field extension of $k$, then we set $X_\ell = X \times_k \ell$.

A variety $X$ over a field $k$ is called rational if $X$ is integral and there is a birational map $X \dashrightarrow \mathbb{P}^n$ for some nonnegative integer $n$; we say that $X$ is geometrically rational if $X_\ell$ is rational for some field extension $\ell$ of $k$.

Recall that the geometric genus $g(C)$ of a geometrically integral curve $C$ over $k$ is defined to be the geometric genus of the unique regular projective geometrically integral model of $C$. If $C$ is itself projective, then this model is the normalisation $\tilde{C}$ of $C$. Note that we have $g(C_{\overline{k}}) \leq g(C)$ with equality if and only if $\tilde{C}$ is smooth (see [Tat52]). In particular, we have $g(C) = 0$ if and only if $\tilde{C}$ is smooth and $C$ is geometrically rational.

Recall that if $C \subset \mathbb{P}^2$ is a curve over $k$ and $S \in C$ a closed point with local ring $\mathcal{O}_{S,C}$, while $C' \subset \mathbb{P}^2$ is a curve that is given in an open neighbourhood $V \subset \mathbb{P}^2$ of $S$ by $h = 0$ for some $h \in \mathcal{O}_{\mathbb{P}^2}(V)$, then the intersection multiplicity $\mu_S(C, C')$ of $C$ and $C'$ at $S$ is the length of the $\mathcal{O}_{S,C}$-module $\mathcal{O}_{S,C}/(h)$, assuming that $S$ does not lie on a common component of $C$ and $C'$. If $S$ is a smooth point of $C$, then the local ring $\mathcal{O}_{S,C}$ is a discrete valuation ring, say with valuation $v_S$, and $\mu_S(C, C')$ equals $v_S(h)$.

We extend the notion of intersection multiplicity, replacing the point $S$ on the curve $C$ by a "branch" of $C$, that is, a point of the normalisation of $C$.

**Definition 5.1.** *Let $C \subset \mathbb{P}^2$ be a curve and let $\vartheta \colon \tilde{C} \to C$ be the normalisation of $C$. Let $T \in \tilde{C}$ be a closed point with local ring $\mathcal{O}_{T,\tilde{C}}$. Let $C' \subset \mathbb{P}^2$ be a curve that is given in an open neighbourhood $V \subset \mathbb{P}^2$ of $\vartheta(T)$ by $h = 0$ for some $h \in \mathcal{O}_{\mathbb{P}^2}(V)$. If the curves $C$ and $C'$ have no irreducible components in common, then the intersection multiplicity $\mu_T(\tilde{C}, C')$ of $\tilde{C}$ and $C'$ at $T$ is defined to be the length of the $\mathcal{O}_{T,\tilde{C}}$-module $\mathcal{O}_{T,\tilde{C}}/(h)$.*

With the same notation as above, the intersection multiplicity $\mu_T(\tilde{C}, C')$ is the same as $\mathrm{ord}_T(h)$ as defined in [Ful98, Section 1.2]. Since $\tilde{C}$ is regular, the local ring $\mathcal{O}_{T,\tilde{C}}$ is a discrete valuation ring, say with valuation $v_T$, and we have $\mu_T(\tilde{C}, C') = v_T(h)$. If $k$ is algebraically closed, then we have $\mu_T(\tilde{C}, C') = \dim_k \mathcal{O}_{T,\tilde{C}}/(h)$.

**Lemma 5.2.** *Let $C, C' \subset \mathbb{P}^2$ be curves and let $\vartheta \colon \tilde{C} \to C$ be the normalisation of $C$. Then for every $S \in C$ we have*
$$\mu_S(C, C') = \sum_{T \mapsto S} \mu_T(\tilde{C}, C') \cdot [k(T) : k(S)],$$
*where the summation runs over all closed points $T \in \tilde{C}$ with $\vartheta(T) = S$ and where $[k(T) : k(S)]$ denotes the degree of the residue field extension.*

*Proof.* This follows immediately from [Ful98, Example 1.2.3]. □

**Definition 5.3.** *Let $C, C' \subset \mathbb{P}^2$ be curves over $k$ that do not have any components in common. Let $\Gamma$ denote either $C$ or its normalisation $\tilde{C}$. Then we define the subset $b(\Gamma, C')$ of $\Gamma(\overline{k})$ as*
$$b(\Gamma, C') = \{T \in \Gamma(\overline{k}) \ : \ \mu_T(\Gamma, C') \text{ is odd}\}.$$

Recall that we are assuming that the characteristic of $k$ is not 2.

**Lemma 5.4.** *Let $D$ be a geometrically integral curve on $X$, let $C = \pi(D)$ be its image under $\pi$, and assume $C$ is not equal to the branch locus $B$. The restriction of $\pi$ to $D$ induces a morphism $\tilde{\pi} \colon \tilde{D} \to \tilde{C}$ between the normalisations of $D$ and $C$. The branch locus of the morphism $\tilde{D}_{\overline{k}} \to \tilde{C}_{\overline{k}}$ induced by $\tilde{\pi}$ is exactly $b(\tilde{C}, B) \subset \tilde{C}(\overline{k})$.*

*Proof.* Without loss of generality, we assume $k = \overline{k}$. Let $\vartheta$ denote the normalisation map $\tilde{C} \to C$. Let $T \in \tilde{C}(k)$ be a point. Since $\tilde{C}$ is regular, the local ring $\mathcal{O}_{T,\tilde{C}}$ is a discrete valuation ring, say with valuation $v_T$. As the characteristic of $k$ is not equal to 2, there is an open neighbourhood $V \subset \mathbb{P}^2$ of $\vartheta(T)$ and an element $h \in \mathcal{O}_{\mathbb{P}^2}(V)$ such that the double cover $\pi^{-1}(V)$ of $V$ is isomorphic to the subvariety of $V \times \mathbb{A}^1(u)$ given by $u^2 = h$. We denote the image of $h$ in the local ring $\mathcal{O}_{T,\tilde{C}}$ and the function field $k(\tilde{C}) = k(C)$ by $h$ as well. The extension $k(C) \subset k(D)$ of function fields is obtained by adjoining a square root $\eta \in k(D)$ of $h$ to $k(C)$. Note that the degree of the restriction of $\pi$ to $D$ is 1 if and only if this extension is trivial, i.e., $h$ is a square in $k(C)$. The intersection $B \cap V$ is given by $h = 0$, so we have $\mu_T(\tilde{C}, B) = v_T(h)$. Suppose $T' \in \tilde{D}(k)$ is a point with $\tilde{\pi}(T') = T$. Since the characteristic of $k$ is not equal to 2, the extension $\mathcal{O}_{T,\tilde{C}} \subset \mathcal{O}_{T',\tilde{D}}$ of discrete valuation rings of $k(C)$ and $k(D) = k(C)(\eta)$, respectively, is ramified if and only if $v_T(h)$ is odd, that is, $T$ is contained in $b(\tilde{C}, B)$, which proves the lemma. □

**Proposition 5.5.** *Let $D$ be a geometrically integral projective curve on $X$, let $C = \pi(D)$ be its image under $\pi$, and assume $g(C) = 0$. Assume also that $C$ is not equal to the branch locus $B$. Let $\tilde{C}$ denote the normalisation of $C$ and set $n = \#b(\tilde{C}, B)$. The following statements hold.*

(1) If $n = 0$, then $\pi$ restricts to a birational morphism $D \to C$ and $g(D) = 0$.
(2) If $n > 0$, then $\pi$ restricts to a double cover $D \to C$ and $g(D) = g(D_{\overline{k}}) = \frac{1}{2}n - 1$.

*Proof.* From $g(C) = 0$, we find that the normalisation $\tilde{C}$ is smooth. Since the characteristic of $k$ is not 2 and for any finite separable field extension $\ell$ of $k$ we have $g(D_\ell) = g(D)$ (see [Tat52, Corollary 2]), we may (and do) replace $k$, without loss of generality, by a quadratic extension $\ell$ for which $\tilde{C}(\ell) \neq \emptyset$. Then $\tilde{C}$ is isomorphic to $\mathbb{P}^1$. Let $\tilde{D}$ denote the normalisation of $D$. The morphism $\pi$ induces a morphism $\tilde{\pi} \colon \tilde{D} \to \tilde{C} \cong \mathbb{P}^1$ of degree at most 2. We claim that $\tilde{D}$ is smooth. Indeed, if $\deg(\tilde{\pi}) = 1$, then this is clear. If $\deg(\tilde{\pi}) = 2$, then because the characteristic of $k$ is not 2, the curve $\tilde{D}$ can be covered by open affine curves that are given by $y^2 = f(x)$ for some polynomial $f \in k[x]$; the regularity of $\tilde{D}$ implies that each polynomial $f$ is separable, which implies that $\tilde{D}$ is smooth. This shows that $g(D) = g(D_{\overline{k}})$, so we may (and do) replace $k$, without loss of generality, by $\overline{k}$.

From $g(C) = 0$, we find that $C$ does not equal the branch locus $B$, so we may apply Lemma 5.4. The Riemann-Hurwitz formula then yields

$$2g(D) - 2 = 2g(\tilde{D}) - 2 = \deg(\tilde{\pi}) \cdot (2g(\tilde{C}) - 2) + n = n - 2\deg(\tilde{\pi}).$$

If $n = 0$, then we find $\deg(\tilde{\pi}) = 1$ and $g(D) = 0$. If $n > 0$, then $\tilde{\pi}$ is not unramified, so $\deg(\tilde{\pi}) = 2$ and we obtain $g(D) = \frac{1}{2}n - 1$. □

**Corollary 5.6.** *Let $C \subset \mathbb{P}^2$ be a geometrically integral projective curve with $g(C) = 0$ that is not equal to the branch locus $B$. Let $\tilde{C}$ denote its normalisation and set $n = \#b(\tilde{C}, B)$. The following statements hold.*

(1) *If $n = 0$, then there exists a field extension $\ell$ of $k$ of degree at most 2 such that the preimage $\pi^{-1}(C_\ell)$ consists of two irreducible components that are birationally equivalent with $C_\ell$. For each such $\ell$ for which $C_\ell$ is rational, the surface $X_\ell$ is unirational.*
(2) *If $n = 0$ and $C$ is rational and there exists a rational point $P \in X(k)$ with $\pi(P) \in C - B$, then the preimage $\pi^{-1}(C)$ consists of two rational components and $X$ is unirational.*
(3) *If $n > 0$, then the preimage $\pi^{-1}(C)$ is geometrically integral and has geometric genus $\frac{1}{2}n - 1$.*
(4) *If $n = 2$ and the preimage $\pi^{-1}(C)$ is rational, then the surface $X$ is unirational.*

*Proof.* As the characteristic of $k$ is not 2, the branch curve $B$ is smooth of genus 3. Set $A = \pi^{-1}(C)$ and $A_{\overline{k}} = A \times_k \overline{k}$. The morphism $A \to C$ induced by $\pi$ has degree 2. From $g(C) = 0$, we find that $C$ is not the branch curve $B$, so $A$ is geometrically reduced. Since $A_{\overline{k}} = A \times_k \overline{k}$ consists of at most two components, there is an extension $\ell$ of $k$ of degree at most 2 such that the components of $A_\ell = A \times_k \ell$ are geometrically irreducible. Let $\ell$ be such an extension and let $D$ be an irreducible component of $A_\ell$.

Suppose $n = 0$. Applying Proposition 5.5 to $D_\ell$ and $C_\ell = \pi(D_\ell)$ shows that the morphism $D_\ell \to C_\ell$ induced by $\pi$ is a birational map. Since $A_\ell \to C_\ell$ has degree 2, there is a unique second component of $A_\ell$, which equals $\iota(D_\ell)$. If, moreover, $C_\ell$ is rational, then Lemma 3.2 implies the unirationality of $X_\ell$, proving statement (1).

Assume $C$ is itself rational and there is a rational point $P \in X(k)$ such that $\pi(P) \in C - B$. Then we have that $P \neq \iota(P)$ and the points $P$ and $\iota(P)$ lie in different components of $A_\ell = D_\ell \cup \iota(D_\ell)$. Since the Galois group $G = G(\ell/k)$ fixes the points $P$ and $\iota(P)$, it follows that $G$ also fixes $D_\ell$ and $\iota(D_\ell)$, so these components are defined over $k$. Then statement (2) follows from (1) taking $\ell = k$.

Suppose $n > 0$. By Proposition 5.5, the morphism $D_{\overline{k}} \to C_{\overline{k}}$ induced by $\pi$ has degree 2, so $D_{\overline{k}}$ is the only component of $A_{\overline{k}}$ and therefore $A$ is geometrically integral. Its genus follows from Proposition 5.5.

Statement (4) follows immediately from Lemma 3.2. □

*Remark* 5.7. Let $D$ be a geometrically integral curve over a field $k$ with $g(D) = 0$. Then there exists a field extension $\ell$ of $k$ of degree at most 2 such that $D_\ell$ is rational. In fact, if $k$ is a finite field, then $D$ is rational over $k$. Therefore, if $k$ is finite in Corollary 5.6, then $C$ is rational; moreover, by case (3) and (4) we conclude that if $n = 2$, then $X$ is unirational over $k$.

*Remark* 5.8. Proposition 5.5 and Corollary 5.6 imply that the geometrically integral projective curves $D \subset X$ with $g(D) = 0$ are exactly the geometrically irreducible components above geometrically integral projective curves $C \subset \mathbb{P}^2$ with $g(C) = 0$ and $\#b(\tilde{C}, B) \in \{0, 2\}$, where $\tilde{C}$ denotes the normalisation of $C$.

5.2. **The pull back of a plane curve in terms of the curve itself.** In the previous subsection, we described the preimage $\pi^{-1}(C) \subset X$ of a curve $C \subset \mathbb{P}^2$, by looking at the intersection points of the branch locus $B$ with the normalisation $\tilde{C}$ of $C$. In this section we will see that we can give an analogous description of $\pi^{-1}(C)$ by looking at the intersection points of $B$ and $C$ itself, if we assume that all singular points of $C$ that lie on $B$ are ordinary singular points.

The following proposition describes the integer $n$ used in Proposition 5.5 in terms of $C$ directly.

**Proposition 5.9.** *Let $C, C' \subset \mathbb{P}^2$ be two projective plane curves with no components in common. Let $\tilde{C}$ be the normalisation of $C$. Assume also that $C'$ is smooth and that all singular points of $C$ that lie on $C'$ are ordinary singularities of $C$. For each point $S \in C(\overline{k})$, let $m_S$ denote the multiplicity of $S$ on $C$. Then we have*

$$\#b(\tilde{C}, C') = \sum_{S \in C(\overline{k}) \cap C'(\overline{k})} c_S(C, C')$$

*with*

$$c_S(C, C') = \begin{cases} m_S & \text{if } m_S \equiv \mu_S(C, C') \pmod{2}, \\ m_S - 1 & \text{if } m_S \not\equiv \mu_S(C, C') \pmod{2}. \end{cases}$$

*Proof.* Let $\vartheta \colon \tilde{C} \to C$ be the normalisation map. Then we may write $b(\tilde{C}, C') = \bigcup_S b_S(\tilde{C}, C')$ with

$$b_S(\tilde{C}, C') = \{T \in \vartheta^{-1}(S) \;:\; \mu_T(\tilde{C}, C') \text{ is odd}\}$$

and where the disjoint union runs over all $S \in C(\overline{k}) \cap C'(\overline{k})$. Suppose $S \in C(\overline{k}) \cap C'(\overline{k})$. Since $C'$ is smooth and the point $S$ is either smooth or an ordinary singularity on $C$, at most one of the $m_S$ points $T \in \vartheta^{-1}(S)$ satisfies $\mu_T(\tilde{C}, C') > 1$. Hence, there is a point $T_0 \in \vartheta^{-1}(S)$ such that for all $T \in \vartheta^{-1}(S)$ with $T \neq T_0$ we have $\mu_T(\tilde{C}, C') = 1$ and thus $T \in b_S(\tilde{C}, C')$. Since we are working over an algebraically closed field, Lemma 5.2 yields $\mu_S(C, C') = \mu_{T_0}(\tilde{C}, C') + m_S - 1$. Hence, we have $T_0 \in b_S(\tilde{C}, C')$ if and only if $m_S$ and $\mu_S(C, C')$ have the same parity. It follows that $\#b_S(\tilde{C}, C') = c_S(C, C')$. The proposition follows. $\square$

We will continue to use the notation $c_S(C, C')$ of Proposition 5.9, which we call the *contribution* of $S$ with respect to $C'$. We set $c_S(C, C') = 0$ for $S \in C(\overline{k})$ with $S \notin C'$.

*Remark* 5.10. Let $C \subset \mathbb{P}^2$ be a geometrically integral projective curve. The points of contribution 0 with respect to $C'$ are the points of $C(\overline{k})$ that are not on $C'$, together with the smooth points $S \in C(\overline{k})$ for which $\mu_S(C, C')$ is even. The points of contribution 1 are the smooth and double points $S$ of $C(\overline{k})$ with $S \in C'$ for which $\mu_S(C, C')$ is odd. The points of type $m > 1$ are the singular points $S$ of $C(\overline{k})$ of multiplicity $m$ or $m + 1$ with $S \in C'$ for which $\mu_S(C, C') \equiv m \pmod{2}$.

Remark 5.8 describes all geometrically integral projective curves $D$ of genus 0 on $X$. In particular, we have $\#b(\tilde{C}, B) \in \{0, 2\}$, where $\tilde{C}$ is the normalisation of the image $C = \pi(D)$. The next proposition can be used to express this condition in terms of $C$ directly, under the assumption that the singular points of $C$ that lie on the branch locus $B$ are ordinary singularities.

**Proposition 5.11.** *Let $C$ and $C'$ be two geometrically integral projective curves in $\mathbb{P}^2$. Let $\tilde{C}$ denote the normalisation of $C$ and let $C^{\mathrm{s}} \subset C(\overline{k})$ denote the set of singular points of $C$. Assume that $C'$ is smooth and that all singular points of $C$ that lie on $C'$ are ordinary. Then the following statements hold.*

*(1) The set $b(\tilde{C}, C')$ is empty if and only if the sets $b(C, C')$ and $C^{\mathrm{s}} \cap C'$ are.*

*(2) We have $\#b(\tilde{C}, C') = 2$ if and only if either*

(a) $b(C, C') = \emptyset$ and there exists a point $S \in C(\overline{k})$ such that $m_S \in \{2, 3\}$ and $C^s \cap C' = \{S\}$, or
(b) there exist two points $S_1, S_2 \in C(\overline{k})$ with $S_1 \neq S_2$ and $b(C, C') = \{S_1, S_2\}$ and $m_{S_1}, m_{S_2} \in \{1, 2\}$ and $C^s \cap C' \subset \{S_1, S_2\}$.

*Proof.* Given that the contributions $c_S(C, C')$ are nonnegative, this follows easily from Proposition 5.9 and Remark 5.10. $\square$

*Remark* 5.12. Suppose $P \in X(k)$ is a rational point that does not lie on the ramification curve, so $\pi(P) \notin B$. Suppose $C$ is a geometrically integral curve of degree $d$ that has a singular point of multiplicity $d - 1$ at $\pi(P)$, and that intersects $B$ with even multiplicity everywhere. Then Proposition 5.11 shows that $b(\tilde{C}, B)$ is empty, so by Corollary 5.6, the pull back $\pi^*(C)$ splits into two components.

If $X$ is general enough, then the Picard group Pic $X$ of $X$ is generated by the canonical divisor $K_X$, and the automorphism group of $X$ acts trivially on Pic $X$, so these two components would be linearly equivalent to the same multiple of $K_X$; as their union is linearly equivalent to $-dK_X$, we find that $d$ is even. Hence, for odd $d$, the answer to the analogous question mentioned below Question 1.4 is negative for $X$ general enough.

It is possible, however, that, even for odd $d$, a variation of this analogous question still has a positive answer. If we forget the del Pezzo surface, and only consider the quartic curve $B \subset \mathbb{P}^2$ with a point $Q \in \mathbb{P}^2$ that does not lie on $B$, we could ask for the existence of a curve of degree $d$ that intersects $B$ with even multiplicity everywhere, and on which $Q$ is a point of multiplicity $d - 1$. The argument above merely shows that if such a curve exists for odd $d$ and $Q$ lifts to a rational point on the del Pezzo surface, then the surface does not have Picard number one.

*Proof of Theorem 1.2.* Assume that the assumptions of statement (1) hold. This implies that $C^s \cap B = \emptyset$ and $b(C, B) = \emptyset$. Therefore, by Proposition 5.11, we have $\#b(\tilde{C}, B) = 0$. Statement (1) follows from applying part (2) of Corollary 5.6.

Assume statement *(2a)* holds. This means that $C^s \cap B = \{Q\}$ and $b(C, B) = \emptyset$. Since $Q$ is a double or triple point of $C$, Proposition 5.11 implies that $\#b(\tilde{C}, B) = 2$. The conclusion of statement (2) follows from applying part (4) of Corollary 5.6 and Remark 5.7.

Assume statement *(2b)* holds. It means that $b(C, B) = \{Q_1, Q_2\}$ and $C^s \cap B \subseteq \{Q_1, Q_2\}$. Since the points $Q_1$ and $Q_2$ are distinct, Proposition 5.11 implies that $\#b(\tilde{C}, B) = 2$. As before, the conclusion of statement (2) follows from part (4) of Corollary 5.6 and Remark 5.7. This concludes the proof of the theorem. $\square$

*Proof of Corollary 1.3.* Set $Q = \pi(P)$. As in Section 4, let $\mathfrak{L}_Q$ denote the line in the dual of $\mathbb{P}^2$ consisting of all lines $L \subset \mathbb{P}^2$ going through $Q$, and note that $\mathfrak{L}_Q$ is isomorphic to $\mathbb{P}^1$. Since $C$ has degree $d$ and $\pi(P)$ is a point of multiplicity $d - 1$, each line in $\mathfrak{L}_Q$ intersects $C$ in a unique $d$-th point, counted with multiplicity. It follows that $C$ is smooth at all points $T \neq Q$. It also follows that the rational map $C \to \mathfrak{L}_Q$ that sends a point $T \in C$ to the line through $T$ and $Q$ is birational, so $C$ is birationally equivalent with $\mathbb{P}^1$. By hypothesis, all intersection points of $B$ and $C$ have even intersection multiplicity.

Assume that $Q$ is not contained in $B$. Since $C$ is smooth away from $Q$, the curve $B$ contains no singular points of $C$. Then $X$ is unirational by part (1) of Theorem 1.2. This proves part (1).

Assume that $Q$ is contained $B$, that $Q$ is an ordinary singularity of $C$, and $d \in \{3, 4\}$. Then $Q$ is a double or a triple point of $C$. Since $Q$ is the only singularity of $C$, the curve $B$ contains no other singular points of $C$. Then $X$ is unirational by part (2) of Theorem 1.2. This proves part (2). $\square$

We now give some examples of curves that satisfy the conditions of Theorem 1.2 or Corollary 1.3.

*Example* 5.13. If $C$ is a bitangent to the branch curve $B$ that is defined over $k$, and $C(k)$ contains a point $Q \notin B$ that lifts to a $k$-rational point on $X$, then Theorem 1.2 implies that $X$ is unirational. We can also prove this directly. Indeed, in this case the pull back $\pi^{-1}(C)$ consists of two exceptional curves that are defined over $k$, so $X$ is not minimal. Blowing down one of these exceptional curves yields a del Pezzo surface $Y$ of degree 3 with a rational point. This implies that $Y$, and therefore also $X$, is unirational.

*Example* 5.14. Suppose the point $P \in X(k)$ is not a generalised Eckardt point and $P$ is not on the ramification curve. Set $Q = \pi(P)$, let $\rho \colon \mathfrak{L}_Q \to X$ be as in Section 4, and set $C = \pi(\rho(\mathfrak{L}_Q))$. Then by Proposition 4.14, the map $\rho(\mathfrak{L}_Q) \to C$ has degree 1, so by Propositions 5.5 and 5.11, the intersection multiplicity of $C$ and the branch curve $B$ is even at all intersection points. Also by Proposition 4.14, the curve $C$ has a point $Q$ off the branch curve $B$ of multiplicity $\deg C - 1$, so the curves of Manin's construction are examples of the curves described in Corollary 1.3. See Remark 4.15.

*Example* 5.15. Consider the surface $X \subset \mathbb{P}(1,1,1,2)$ over $\mathbb{F}_3$, defined by the equation
$$w^2 = x^4 + y^4 + z^4.$$
The surface $X$ is a del Pezzo surface of degree 2. All its rational points either are on the ramification curve, or they are generalised Eckardt points. In fact, the surface $X$ has 154 rational points over $\mathbb{F}_9$, with 28 of those lying on the ramification locus. The remaining 126 are generalised Eckardt points, which is also the maximum number of generalised Eckardt points a del Pezzo surface of degree two can have (see [STVA14, before Example 7]). It follows that Manin's method does not apply to this surface. Let $P$ be the point $(0:0:1:1)$ on $X$. Then $P$ is a generalised Eckardt point and its image $Q = \pi(P) = (0:0:1) \in \mathbb{P}^2$ does not lie on the branch locus $B$, which is given by $x^4 + y^4 + z^4 = 0$. Consider the curve $C \subset \mathbb{P}^2$ given by $x^3y + xy^3 = z(x+y)^2(y-x)$. The curve $C$ is a geometrically integral quartic plane curve that has a triple point at $Q$ and that intersects $B$ with even multiplicity everywhere. Therefore, by case (1) of Corollary 1.3, the surface $X$ is unirational.

Of course, unirationality of $X$ was already known: it follows for instance from Lemma 20 in [STVA14] (cf. Example 5.17 below). It is nice to see, though, that, even though Manin's construction and the generalisation in [STVA14] do not produce a curve in $\mathbb{P}^2$ of some degree $d$ with a point of multiplicity $d-1$ at $Q$, and even intersection multiplicity with $B$ everywhere, such curves do still exist, and then case (1) of Corollary 1.3 implies unirationality of $X$. This gives a positive answer to the question below Question 1.4 for $d = 4$ and this particular surface $X$ and this generalised Eckardt point $P$.

One might ask whether there are curves of lower degree satisfying the hypotheses of case (1) of Corollary 1.3. Indeed, there are conics that do, for example the one given by $y^2 = xz$. An exhaustive computer search, based on Corollary 5.6, part (2), and Corollary 6.2, shows that there are no cubic curves with a double point at $Q$ satisfying the hypotheses of Corollary 1.3 and its case (1).

*Example* 5.16. Let $X_1, X_2, X_3$ be the three del Pezzo surfaces defined as in Section 3 and let $B_i$ be their branch locus, for $i = 1, 2, 3$. For $i = 1, 2, 3$, all rational points of the surface $X_i$ lie on the ramification locus. Consider the rational points $P_1 = (0:0:1:0) \in X_1$, $P_2 = (0:1:1:0) \in X_2$, and $P_3 = (1:1:1:0) \in X_3$, and set $Q_i = \pi(P_i)$. Clearly, we have $Q_i \in B_i$. Set $d_1 = d_2 = 4$ and $d_3 = 3$. Let $C_i \subset \mathbb{P}^2$ be the projective plane curve of degree $d_i$ given by the polynomial $h_i$ defined as in Remark 3.4. The curve $C_i$ is geometrically irreducible and it has an ordinary singular point at $Q_i$ of multiplicity $d_i - 1$. Given that the curve $C_i$ pulls back to the geometrically irreducible rational curve $A_i$ of Remark 3.4, we find from Corollary 5.6 and Proposition 5.11 that $C_i$ intersects $B_i$ with even multiplicity everywhere.

Of course, one could also check directly that $C_i$ intersects $B_i$ with even multiplicity everywhere. Then Corollary 1.3 and Remark 5.7 give an alternative proof that the surface $X_i$ is unirational (cf. Remark 3.4). There is a quartic alternative for $C_3$ as well. The curve $C_3' \subset \mathbb{P}^2$ given by the vanishing of
$$\begin{aligned}h_3' =\,& \gamma^2 x^4 + x^3y + \gamma x^2y^2 + \gamma^3 xy^3 - y^4 + x^3z + \gamma x^2yz + xy^2z \\ & - \gamma y^3z + \gamma x^2z^2 + xyz^2 + \gamma^3 y^2z^2 + \gamma^3 xz^3 - \gamma yz^3 - z^4\end{aligned}$$
is geometrically integral, has an ordinary triple point at $(-1:1:1)$, and intersects $B$ with even multiplicity everywhere.

*Example* 5.17. Let $k$ be a field with characteristic different from 2. Let $a_1, \ldots, a_6 \in k$ be such that the variety $X$ in the weighted projective space $\mathbb{P} = \mathbb{P}(1,1,1,2)$ defined by

$$w^2 = a_1^2 x^4 + a_2^2 y^4 + a_3^2 z^4 + a_4 x^2 y^2 + a_5 x^2 z^2 + a_6 y^2 z^2$$

is a del Pezzo surface of degree 2. This is the surface of Lemma 20 in [STVA14], where it is noted that the surface in $\mathbb{P}$ given by $w = a_1 x^2 + a_2 y^2 + a_3 z^2$ intersects the surface $X$ in a curve $D$, which the anti-canonical map $\pi \colon X \to \mathbb{P}^2$ sends isomorphically to the plane quartic curve $C \subset \mathbb{P}^2$ given by

$$(a_4 - 2a_1 a_2) x^2 y^2 + (a_5 - 2a_1 a_3) x^2 z^2 + (a_6 - 2a_2 a_3) y^2 z^2 = 0.$$

They also note that this curve $C$ is birationally equivalent to a conic under the standard Cremona transformation, so $C$ and $D$ are rational over an extension of $k$ of degree at most 2. If they are rational over $k$, then $X$ is unirational.

Indeed, one checks that the curve $C$ satisfies the conditions of part (1) of Corollary 5.6, and if $C$ is rational over $k$, then it also satisfies the conditions of part (1) of Theorem 1.2, where one can take $P$ to be any of the points on $X$ above any of the singular points $(0:0:1)$, $(0:1:0)$, and $(1:0:0)$ of $C$.

## 6. Finding appropriate curves

In this section, we assume that the characteristic of $k$ is not 2, and we give sufficient easily-verifiable conditions for a curve $C$ to satisfy the hypotheses of Corollary 1.3. This is also how we found the three curves, $C_1, C_2$, and $C_3$ of Remark 3.4, whose existence implies unirationality of the three difficult surfaces $X_1$, $X_2$, $X_3$ (see Example 5.16 and Remark 6.7).

Let $X \subset \mathbb{P}(1,1,1,2)$ be a del Pezzo surface of degree 2, given by $w^2 = g$ with $g \in k[x,y,z]$ homogeneous of degree 4. Let $B \subset \mathbb{P}^2(x,y,z)$ be the branch curve of the projection $\pi \colon X \to \mathbb{P}^2$. Then $B$ is given by $g = 0$. Let $P \in X(k)$ be a rational point and set $Q = \pi(P)$. Without loss of generality, we assume $Q = (0:0:1)$. Let $C \subset \mathbb{P}^2$ be a geometrically irreducible curve of degree $d \geq 2$ on which $Q$ is a point of multiplicity $d - 1$.

There are coprime homogeneous polynomials $f_{d-1}, f_d \in k[x,y]$ of degree $d-1$ and $d$, respectively, such that $C$ is given by $z f_{d-1} = f_d$. The projection away from $Q$ induces a birational map from $C$ to the family $\mathfrak{L}_Q$ of lines in $\mathbb{P}^2$ through $Q$. Its inverse is a morphism $\vartheta$ that sends a line $L \in \mathfrak{L}_Q$ to the $d$-th intersection point of $L$ with $C$. If we identify $\mathfrak{L}_Q$ with $\mathbb{P}^1$, where $(s:t) \in \mathbb{P}^1$ corresponds to the line given by $sy = tx$, then $\vartheta \colon \mathbb{P}^1 \to C$ sends $(s:t)$ to

$$(s f_{d-1}(s,t) : t f_{d-1}(s,t) : f_d(s,t)).$$

The curve $C$ has no singularities outside $Q$, and we may identify $\vartheta \colon \mathbb{P}^1 \to C$ with the normalisation of $C$. The points on $\mathbb{P}^1$ above the point $Q$ are exactly the points where $f_{d-1}(s,t)$ vanishes. The curve $C$ has an ordinary singularity at $Q$ if and only if $d > 2$ and $f_{d-1}(s,t)$ vanishes at $d-1$ distinct $\overline{k}$-points of $\mathbb{P}^1(s,t)$.

The pull back $\pi^*(C)$ is birationally equivalent with the curve given by $w^2 = G$ in the weighted projective space $\mathbb{P}(1,1,2d)$ with coordinates $s, t, w$, and with

$$G = g\big(s f_{d-1}(s,t), t f_{d-1}(s,t), f_d(s,t)\big) \in k[s,t].$$

**Proposition 6.1.** *For any point $T \in \mathbb{P}^1(\overline{k})$, the intersection multiplicity $\mu_T(\mathbb{P}^1, B)$ equals the order of vanishing of $G$ at $T$.*

*Proof.* Since $C$ either has degree 2 or it is singular, it is not equal to $B$. As $C$ is irreducible, it has no irreducible components in common with $B$. By symmetry between $s$ and $t$, we may assume $T = (\alpha : 1)$ for some $\alpha \in \overline{k}$. Then the local ring $\mathcal{O}_{T,\mathbb{P}^1}$ is isomorphic the the localisation of $\overline{k}[s]$ at the maximal ideal $(s - \alpha)$. Let $\ell \in k[x,y,z]$ be a linear form that does not vanish at $\vartheta(T)$. Then locally around $\vartheta(T) \in \mathbb{P}^2$, the curve $B$ is given by the vanishing of the element $g/\ell^4$, whose image in $\mathcal{O}_{T,\mathbb{P}^1}$ is $G(s,1)/L(s,1)^4$ with $L(s,t) = \ell\big(s f_{d-1}(s,t), t f_{d-1}(s,t), f_d(s,t)\big)$. Since $L(s,1)$ does not vanish at $\alpha$, we find that $\mu_T(\mathbb{P}^1, B)$ equals the order of vanishing of $G(s,1)$ at $\alpha$, which equals the order of vanishing of $G$ at $T$. $\square$

**Corollary 6.2.** *We have $b(\mathbb{P}^1, B) = \emptyset$ if and only if $G$ is a square in $\overline{k}[s,t]$.*

*Proof.* By Proposition 6.1, we have $b(\mathbb{P}^1, B) = \emptyset$ if and only if the order of vanishing of $G$ is even at every point $T \in \mathbb{P}^1(\overline{k})$. This is equivalent with $G$ being a square in $\overline{k}[s,t]$. □

If $B$ does not contain the unique singular point $Q$ of $C$, then $\vartheta$ induces a bijection $b(\mathbb{P}^1, B) \to b(C, B)$, so in this case we also have $b(C, B) = \emptyset$ if and only if $G$ is a square in $\overline{k}[s,t]$. The following proposition gives an analog of this statement when $Q$ is contained in $B$.

**Proposition 6.3.** *Suppose that $Q$ is contained in $B$. Then $H = G/f_{d-1}(s,t)$ is contained in $k[s,t]$. Suppose, furthermore, that the tangent line to $B$ at $Q$ is given by $h = 0$ with $h \in k[x,y]$, and that $Q$ is an ordinary singular point on the curve $C$. Then the following statements hold.*
  (1) *Suppose $d$ is odd. Then the set $b(C, B)$ is empty if and only if $H$ is a square in $\overline{k}[s,t]$.*
  (2) *Suppose $d$ is even. If $h$ divides $f_{d-1}$, then $H/h(s,t)$ is contained in $k[s,t]$. The set $b(C,B)$ is empty if and only if $h$ divides $f_{d-1}$ and $H/h(s,t)$ is a square in $\overline{k}[s,t]$.*

*Proof.* Write $g = \sum_{i=0}^{4} g_i z^{4-i}$, where $g_i \in k[x,y]$ is homogeneous of degree $i$ for all $0 \leq i \leq 4$. If $g(Q) = g_0$ vanishes, then each monomial of $g$ is divisible by $x$ or $y$, which implies that $G$ is divisible by $f_{d-1}$, which in turn shows $H \in k[s,t]$. Suppose that all hypotheses hold. By $g(Q) = 0$ we find $g_0 = 0$. The tangent line to $B$ at $Q$ is given by $g_1 = 0$, so $h$ is a scalar multiple of $g_1$. Note that all statements are invariant under the action of $\mathrm{GL}_2(k)$ on $\mathbb{P}^1$ and $\mathbb{P}^2$ given on their respective homogeneous coordinate rings $k[s,t]$ and $k[x,y,z]$ by $\gamma(s) = as+bt, \gamma(t) = cs+dt$ and $\gamma(x) = ax+by, \gamma(y) = cx+dy, \gamma(z) = z$ for
$$\gamma = \begin{pmatrix} a & b \\ c & d \end{pmatrix}.$$
After applying an appropriate element $\gamma \in \mathrm{GL}_2(k)$ and rescaling $h$, we assume, without loss of generality, that $h = g_1 = y$.

If $y$ divides $f_{d-1}$, then $t$ divides $f_{d-1}(s,t)$; since all monomials in $g$ besides $y$ are divisible by $x^2$, $xy$, or $y^2$, it follows that in this case $G$ is divisible by $tf_{d-1}(s,t)$, so $H/t$ is contained in $k[s,t]$. This does not depend on $d$ being even.

In the open neighbourhood of $Q$ given by $z \neq 0$, the curve $B$ is given by the vanishing of $g/z^4 = g(x/z, y/z, 1) = \sum_i g_i(x/z, y/z)$. The maximal ideal $\mathfrak{m}$ of the local ring $\mathcal{O}_{Q,C}$ is generated by $x/z$ and $y/z$, so the image of $g/z^4$ in $\mathcal{O}_{Q,C}/\mathfrak{m}^2$ is $g_1(x/z, y/z) = y/z$. Let $T \in \mathbb{P}^1(\overline{k})$ be a point with $\vartheta(T) = Q$, and let $\mathfrak{n}$ be the maximal ideal of the local ring $\mathcal{O}_{T,\mathbb{P}^1}$. Then the image of $g$ in $\mathcal{O}_{T,\mathbb{P}^1}/\mathfrak{n}^2$ equals the image of $y/z$, which is $tf_{d-1}(s,t)/f_d(s,t)$. The point $T$ corresponds to a linear factor of $f_{d-1}(s,t)$. Since $f_d(s,t)$ does not vanish at $T$, we find that the valuation $v_T(g)$ of $g$ in $\mathcal{O}_{T,\mathbb{P}^1}$ is at least 2 if $t$ vanishes at $T$, that is, $\mu_T(\mathbb{P}^1, B) \geq 2$ if $T = (1:0)$. We have $\mu_T(\mathbb{P}^1, B) = v(g) = 1$ if $T \neq (1:0)$. From Lemma 5.2 we conclude

(5) $$\mu_Q(C, B) = \begin{cases} d-2 + \mu_{(1:0)}(\mathbb{P}^1, B) & \text{if } y \text{ divides } f_{d-1}, \\ d-1 & \text{otherwise.} \end{cases}$$

We now consider the two cases.
  (1) Suppose $d$ is odd. From (5) it follows that $\mu_Q(C, B)$ is even if and only if either $y$ divides $f_{d-1}$ and $\mu_{(1:0)}(\mathbb{P}^1, B)$ is odd, or $y$ does not divide $f_{d-1}$. This happens if and only if $\mu_T(\mathbb{P}^1, B)$ is odd for all $T \in \mathbb{P}^1$ at which $f_{d-1}(s,t)$ vanishes. For all other points $R \in C$ with $R \neq Q$, the multiplicity $\mu_R(C, B)$ is even if and only if $\mu_{\vartheta^{-1}(R)}(\mathbb{P}^1, B)$ is even. From Proposition 6.1, we conclude that $b(C, B)$ is empty if and only if the order of vanishing of $G$ is odd at all points $T \in \mathbb{P}^1$ at which $f_{d-1}(s,t)$ vanishes and even at all other points. This is equivalent with $H$ being a square in $\overline{k}[s,t]$.
  (2) Suppose $d$ is even. From (5) it follows that $\mu_Q(C, B)$ is even if and only if $y$ divides $f_{d-1}$ and $\mu_{(1:0)}(\mathbb{P}^1, B)$ is even. As in the case for odd $d$, this implies that $b(C, B)$ is empty if and only if the order of vanishing of $G$ is odd at all points $T \neq (1:0)$ at which $f_{d-1}(s,t)$ vanishes, and even at all other points, including $(1:0)$. Since the order of vanishing of $tf_{d-1}(s,t)$ at $(1:0)$ is 2, this is equivalent to $G/(tf_{d-1}(s,t)) = H/t$ being a square in $\overline{k}[s,t]$.

This finishes the proof. □

We have already seen that the pull back $\pi^*(C)$ is birationally equivalent with the curve given by $w^2 = G$ in $\mathbb{P}(1, 1, 2d)$. This curve splits into two $k$-rational components if and only if $G$ is a square in $k[s, t]$. If $Q$ is an ordinary singular point of $C$ that lies on $B$, then this never happens. However, the curve $\pi^*(C)$ may itself be $k$-rational, in which case $G$ factors as a square times a quadric.

We will now focus on the case $d = 4$, so $Q$ is a triple point. The following corollary says that if $Q$ is an ordinary triple point, then we do not need to factorise $G$, as we know exactly which part should be the square, and which the quadric.

**Corollary 6.4.** *Suppose that $Q$ is an ordinary singular point of $C$ that lies on $B$. If the pull back $\pi^*(C) \subset X$ is $k$-rational, then we have $d \leq 4$.*

*Moreover, suppose $d = 4$, and let the tangent line to $B$ at $Q$ be given by $h = 0$ with $h \in k[x, y]$. Then the pull back $\pi^*(C) \subset X$ is $k$-rational if and only if there is a constant $c \in k^*$ such that the following statements hold:*

  *(1) the polynomial $h$ divides $f_3$;*
  *(2) the polynomial $cH(s, t)/h(s, t)$ is a square in $k[s, t]$;*
  *(3) the conic given by $cw^2 = f_3(s, t)/h(s, t)$ in $\mathbb{P}^2(s, t, w)$ is $k$-rational.*

*Proof.* Suppose $\pi^*(C)$ is $k$-rational. Then $\pi^*(C)$ is geometrically integral and has genus $g(\pi^*(C)) = 0$. From Proposition 5.5 we obtain $b(\mathbb{P}^1, B) = 2$. From Proposition 5.9 we conclude that the contribution $c_Q(C, B)$ is at most 2. Moreover, this proposition also gives $c_Q(C, B) \geq d - 2$ with equality if and only if $\mu_Q(C, B)$ is even. We conclude $d \leq 4$.

Suppose $d = 4$. Then we have equality $c_Q(C, B) = 2 = \#b(\mathbb{P}^1, B)$, so $\mu_Q(C, B)$ is even, and we find that $b(C, B)$ is empty. From Proposition 6.3 we find that $h$ divides $f_3$, and $m = H(s, t)/h(s, t)$ is a square in $\overline{k}[s, t]$. Let $c$ be the main coefficient of $m(s, 1)$. Then $cm$ is a square in $k[s, t]$. Therefore, the $k$-rational curve given by $w^2 = G$ with

$$(6) \qquad G = cm \cdot h^2(s, t) \cdot c^{-1} f_3(s, t)/h(s, t)$$

in $\mathbb{P}(1, 1, 2d)$ is birationally equivalent with the conic given by $cw^2 = f_3(s, t)/h(s, t)$ in $\mathbb{P}^2(s, t, w)$, which is therefore also $k$-rational.

Conversely, if there is a constant $c$ such that $cH(s, t)/h(s, t)$ is a square in $k[s, t]$, then it follows from (6) that the conic given by $cw^2 = f_3(s, t)/h(s, t)$ in $\mathbb{P}^2(s, t, w)$ is birationally equivalent with the curve in $\mathbb{P}(1, 1, 2d)$ given by $w^2 = G$, which is birationally equivalent with $\pi^*(C)$. Hence, if this conic is $k$-rational, then so is $\pi^*(C)$. □

*Remark* 6.5. Corollary 6.4 helps us in finding all curves $C$ of degree $d = 4$ that satisfy the conditions of case (2) of Corollary 1.3 with $\ell = k$. More explicitly, after a linear transformation of $\mathbb{P}^2$, we may assume that $Q = (0 : 0 : 1)$, and the tangent line to $B$ at $Q$ is given by $y = 0$. Then we claim that every curve $C$ of degree $d = 4$ that satisfies the conditions of case (2) of Corollary 1.3 with $\ell = k$ is given by

$$yz\phi_2 = x^4 + y\phi_3$$

for some homogeneous $\phi_2, \phi_3 \in k[x, y]$ of degree 2 and 3, respectively, with $\phi_2$ squarefree and not divisible by $y$. Indeed, we find that $f_3$ is divisible by $y$, so there is a $\phi_2 \in k[x, y]$ such that $f_3 = y\phi_2$; since $C$ is irreducible, the polynomial $f_4$ is not divisible by $y$, so the coefficient of $x^4$ in $f_4$ is nonzero, and after scaling $\phi_2, f_3$, and $f_4$, we may assume that there exists a $\phi_3 \in k[x, y]$ such that $f_4 = x^4 + y\phi_3$. Moreover, $Q$ is an ordinary singularity if and only if $\phi_2$ is squarefree and not divisible by $y$.

Hence, to find all such curves $C$, we are looking for all pairs $(\phi_2, \phi_3)$ with $\phi_i \in k[x, y]$ homogeneous of degree $i$, such that

  (1) the polynomial $\phi_2$ is squarefree and $y$ does not divide $\phi_2$,
  (2) the curve given by $yz\phi_2 = x^4 + y\phi_3$ is geometrically integral,
  (3) there is a constant $c \in k^*$ such that polynomial $c \cdot G(s, t)/(t^2\phi_2(s, t))$ with

$$G = g\big(st\phi_2(s, t),\ t^2\phi_2(s, t),\ s^4 + t\phi_3(s, t)\big)$$

is a square,
(4) the conic given by $cw^2 = \phi_2(s,t)$ in $\mathbb{P}^2(s,t,w)$, with $c$ as in (3), is $k$-rational.

Note for (3) that, because the characteristic is not 2, a homogeneous polynomial $H \in k[s,t]$ of even degree is a square in $\overline{k}[s,t]$ if and only if there is a constant $c \in k^*$ such that $cH$ is a square in $k[s,t]$, which happens if and only if $\gamma^{-1}H(s,1)$ is a square in $k[s]$, where $\gamma$ is the main coefficient of $H(s,1)$. This follows from the fact that a monic polynomial in $k[s]$ is a square in $k[s]$ if and only if it is a square in $\overline{k}[s]$. Moreover, the $c \in k^*$ for which $cH$ is a square, form a coset in $k^*/k^{*2}$, so whether or not (4) holds does not depends on the choice of $c$.

Question 1.4 for $d = 4$ can be rephrased using Remark 6.5. It is equivalent to the following question.

**Question 6.6.** *Let $k$ be a field of characteristic not equal to 2, and $g \in k[x,y,z]$ a homogeneous polynomial of degree 4 such that the curve $B \subset \mathbb{P}^2(x,y,z)$ given by $g = 0$ is smooth, it contains the point $Q = (0:0:1)$, and the tangent line to $B$ at $Q$ is given by $y = 0$. Do there exist homogeneous polynomials $\phi_2, \phi_3 \in k[x,y]$ of degree 2 and 3, such that conditions (1)–(3) of Remark 6.5 are satisfied?*

*Remark* 6.7. If $k$ is a (small) finite field, then we can list all pairs $(\phi_2, \phi_3)$ with $\phi_i \in k[x,y]$ homogeneous of degree $i$, and check for each whether the conditions (1)–(4) of Remark 6.5 are satisfied. In fact, condition (4) is automatically satisfied over finite fields. Indeed, this is how we found the curves $C_1, C_2$ given in Remark 3.4, whose existence implies unirationality of the three difficult surfaces $X_1, X_2$ (see Example 5.16). Finding the rational cubic curve $C_3$ on $X_3$, as given in Remark 3.4, was easier, based on part (1) of Proposition 6.3.

*Remark* 6.8. For any integer $i$, let $k[x,y]_i$ denote the $(i+1)$-dimensional space of homogeneous polynomials of degree $i$. In general, over any field, we can describe the set of pairs $(\phi_2, \phi_3) \in k[x,y]_2 \times k[x,y]_3$ satisfying condition (3) of Remark 6.5 as follows.

Identify $k[x,y]_2 \times k[x,y]_3$ with the affine space $\mathbb{A}^7$ and let $R$ denote the coordinate ring of $\mathbb{A}^7$, that is, $R$ is the polynomial ring in the $3 + 4 = 7$ coefficients of $\phi_2$ and $\phi_3$. Let $Z \subset \mathbb{A}^7$ be the locus of all $(\phi_2, \phi_3)$ that satisfy condition (3).

For *generic* $\phi_2, \phi_3$, that is, with the variables of $R$ as coefficients, the coefficients of the polynomial
$$G' = G(s,t)/(t^2 \phi_2(s,t))$$
of condition (3) of Remark 6.5 lie in $R$. For general enough $g$, the coefficient $c \in R$ of $s^{12}$ in the polynomial $G' \in R[s,t]$ is nonzero. On the open set $U$ of $\mathbb{A}^7$ given by $c \neq 0$, we may complete $G'(s,1)$ to a square in the sense that there are polynomials $G_1, G_2 \in R[c^{-1}][s]$ with $G_1$ monic of degree 6 in $s$ and $G_2$ of degree at most 5 in $s$ such that $G'(s,1) = cG_1^2 - G_2$. The vanishing of the six coefficients in $R$ of $G_2$ determines the locus $Z \cap U$ inside $U$ of all pairs $(\phi_2, \phi_3)$ at which $cG'$ is a square. Note that we have $c = G'(1,0)$. For each point $(s_0 : t_0) \in \mathbb{P}^1$, we can use an automorphism of $\mathbb{P}^1$ that sends $(s_0 : t_0)$ to $(1 : 0)$, to similarly describe the intersection of $Z$ with the open subset of $\mathbb{A}^7$ where $G'(s_0, t_0)$ is nonzero; it is also given by the vanishing of six polynomials in $R$. We can cover $\mathbb{A}^7$ with open subsets of this form, thus describing $Z$ completely.

A naive dimension count suggests that the locus $Z$ has dimension $7 - 6 = 1$. This is consistent with the following, similarly naive, dimension count. The family of quartic curves in $\mathbb{P}^2$ is 14-dimensional, as it is the projective space $\mathbb{P}(k[x,y,z]_4)$, where $k[x,y,z]_4$ is the 15-dimensional vector space of polynomials of degree 4. The codimension of the subset of those curves having a triple point at $Q$ is 6, and demanding that the intersection multiplicity $\mu_Q(C, B)$ is at least 4 cuts down another dimension. Since $B$ is also a quartic curve, by Bezout's theorem it follows that $B$ and $C$ have 16 intersection points, counted with multiplicity. Hence, generically, the curves in the remaining 7-dimensional family intersect $B$, besides in $Q$, in $16 - 4 = 12$ more points. One might expect the subfamily of those curves where this degenerates to six points with multiplicity 2 to have codimension 6, in which case this would leave a 1-dimensional family of quartic curves with a triple point at $Q \in B$ and intersecting $B$ with even multiplicity everywhere.

However, the locus $Z$ also contains some degenerate components that we are not interested in. For example, the locus of all $(0, \phi_3)$ for which $f_4 = x^4 + y\phi_3$ is a square is contained in $Z$ and has dimension 2. Also, for any smooth conic $\Gamma$ that contains $Q$, that has its tangent line at $Q$ given by $y = 0$, and that has even intersection multiplicity with $B$ everywhere, we get a 1-dimensional subset of $Z$ consisting of pairs $(\phi_2, \phi_3)$ that correspond with the union of $\Gamma$ with any double line through $Q$ (these lines are parametrised by $\mathbb{P}^1$). Note that in all these degenerate cases the curve $C$ is reducible. Another degenerate case is the limit of Manin's construction. By Remark 4.11, this limit curve is the non-reduced curve $\pi_*\pi^*(2L) = 4L$, where $L$ is the tangent line to $B$ at $Q$, given by $y = 0$. Hence, this quartic curve is given by $y^4 = 0$, which does not correspond to a point on the affine set $Z$, as the coefficient of $x^4$ is zero.

Let $Z_0$ denote the affine subset of $Z$ corresponding to curves $C$ that are geometrically integral and on which $Q$ is an ordinary triple point. Then Questions 1.4 (for $d = 4$) and 6.6 can be rephrased by asking whether the subset $Z_0$ contains a $k$-rational point.

*Example* 6.9. Let $B \subset \mathbb{P}^2$ be the smooth curve given by
$$y^4 - x^4 - x^3y - xy^3 + y^3z + yz^3 = 0$$
over $k = \mathbb{F}_3$, and let $Q$ be the point $(0 : 0 : 1) \in B(k)$. The tangent line to $B$ at $Q$ is given by $y = 0$. We ran through all the homogeneous polynomials in $x, y$ of degree 2 and 3 over $k$ and we did not find any pair $(\phi_2, \phi_3)$ of polynomials satisfying conditions (1)–(3) of Remark 6.5; we did find pairs satisfying only conditions (2) and (3). This means that Questions 1.4 (for $d = 4$) and 6.6 have negative answer in this specific case. It could, however, still be true that the answer is positive for $X$ and $P$ general enough.

Notice that the curve $B$ is isomorphic to the Fermat curve $x^4 + y^4 + z^4 = 0$ via the following linear change of variables
$$(x : y : z) \mapsto (x - z : x + y + z : z).$$